\documentclass[a4paper, oneside, 11pt]{article}
\usepackage[english]{babel}
\usepackage{amsfonts}
\usepackage{amssymb}
\usepackage{graphics}
\usepackage{amsrefs}
\usepackage{latexsym}

\RequirePackage[OT1]{fontenc}
\RequirePackage{amsthm,amsmath}
\RequirePackage[numbers]{natbib}
\RequirePackage[colorlinks,citecolor=blue,urlcolor=blue]{hyperref}
\usepackage{breakurl}
\usepackage{amssymb}
\usepackage[latin1]{inputenc}
\theoremstyle{plain}
\newtheorem{theorem}{Theorem}[]
\newtheorem{proposition}{Proposition}[]

\theoremstyle{remark}
\newtheorem{remark}{Remark}[]
\newtheorem{keywords}{Some key words}
\theoremstyle{definition}
\newtheorem{example}{Example}[]

\begin{document}

%\jname{Biometrika}
%% The year, volume, and number are determined on publication
%\jyear{201X}
%\jvol{XX}
%\jnum{1}
%% The \doi{...} and \accessdate commands are used by the production team
%\doi{10.1093/biomet/asm023}
%\accessdate{Advance Access publication on 31 July 2012}
%\copyrightinfo{\Copyright\ 2012 Biometrika Trust\goodbreak {Em Printed in Great Britain}}

%% These dates are usually set by the production team
%\received{April 201X}
%\revised{September 201X}

%% The left and right page headers are defined here:
%\markboth{A. Cerquetti}{BNP estimation of Tsallis diversity}

%% Here are the title, author names and addresses
%\title{Bayesian nonparametric estimation of Tsallis diversity indices under Gnedin--Pitman priors}

%\author{ANNALISA CERQUETTI}
%\affil{Department of Methods and Models for Economics, Territory and Finance,\\ Sapienza University of Rome, Via del Castro Laurenziano, 9,\\ 00161 Rome, Italy \email{annalisa.cerquetti@gmail.com}}

%\author{\and D. M. TITTERINGTON}
%\affil{Department of Statistics, University of Glasgow, Glasgow G12 8QQ, U.K. Email{mike@stats.gla.ac.uk}}

%\maketitle 

\title{{\bf{\Large{Bayesian nonparametric estimation of Tsallis diversity indices under Gnedin--Pitman priors}}}\footnote{{\it AMS (2000) subject classification}. Primary: 60G58. Secondary: 60G09.}}
\author{\textsc {Annalisa Cerquetti}\footnote{Corresponding author, SAPIENZA University of Rome, Via del Castro Laurenziano, 9, 00161 Rome, Italy. E-mail: {\tt annalisa.cerquetti@gmail.com}}\\
\it{\small Department of Methods and Models for Economics, Territory and Finance}\\
  \it{\small Sapienza University of Rome, Italy }}
\date{\today}
\maketitle{}

\begin{abstract}
Tsallis entropy is a generalized diversity index first derived in Patil and Taillie (1982) and then rediscovered in community ecology by Keylock (2005).
%as a concave diversity measure addressing the self-similar and multifractal nature of the distribution and abundance of species. 
Bayesian nonparametric estimation of Shannon entropy and Simpson's diversity under uniform and symmetric Dirichlet priors has been already advocated as an alternative to maximum likelihood estimation based on frequency counts, which is negatively biased in the undersampled regime. %Two recent results provides Bayesian estimation of Shannon and Simpson indexes  under two-parameter Poisson-Dirichlet priors. 
Here we present a fully general Bayesian nonparametric estimation of the whole class of  Tsallis diversity indices under Gnedin-Pitman priors, a large family of random discrete distributions 
%generalizing the Poisson--Dirichlet distribution and 
recently deeply investigated in posterior predictive species richness and discovery probability estimation. We provide both prior and posterior analysis. The results, illustrated through examples and an application to a real dataset, show the procedure is easily implementable, flexible and overcomes limitations of previous frequentist and Bayesian solutions.
%Our results complete the large set of results in Bayesian nonparametrics  for species sampling problems restricted until now to species richness and discovery probability estimation.
\end{abstract}

\begin{keywords}
Bayesian nonparametrics; Diversity; Entropy; Gnedin-Pitman priors; Shannon entropy; Simpson's diversity; Species sampling; Tsallis entropy.
\end{keywords}

\section{Introduction}
\subsection{A generalized diversity index}
%Species diversity is made of two components:  In other words a population of species is the more {\it diverse} the more is {\it rich} and {\it uniform}. 
Diversity is both a goal and an indicator of ecosystems health and function. The measurement of diversity of populations when individuals are classified into groups has a long history, dating back to Simpson's (1949) and Fisher's (1943) seminal papers. Since then the ecological literature has produced a variety of indices to measure both species {richness}, the number of different species belonging to a population, and species {evenness}, the distance of the actual relative abundances from a situation of uniform distribution of the population into different species. 
%Being both those quantities commonly adopted both as goals and indicators of ecosystem health and function (... ).
In 1982 Patil and Taillie generalize Shannon entropy (Shannon, 1948) and Simpson's index, by far the most widely used measures of biological diversity, identifying a generalized diversity measure as the mean value
%for $m\neq 1$, 
%or
%$$
%R(P_i)= -\log(P_i) 
%$$
%for  $m=1$ 
%thus defining 
\begin{equation}
\label{tsallis}
H_m(P)= \frac{1}{m-1} \left(1- \sum_i P_i^m\right),
\end{equation}  
%the solving the proportionality equation
%by deriving 
%a dichotomous rarity of species $R(P_i)$ 
%solution to the {\it proportionality equation}
%$
%R(qP_i)-R(P_i)=W(P_i)[R(q)-R(1)],
%$
for $P=(P_i)_{i \geq 1}$ a population of relative abundances and $m >0$ a parameter specifying 
%$m > 0$ and $m \neq 1$.
%$R(P_i)$a rarity function, $q, R, W \in (0,1]$ and $W(P_i)$ a {\it deflation factor}.
% and $R(P_i)$ a dichotomous rarity index, i.e. one depending only on the relative abundance of the different species. 
%Apart from degenerate cases, they show that, for $W(P_i)=P_i^{m-1}$, possible solutions for correspond to
%$
%R(P_i)= (1-P_i^{m-1})/(m-1)  
%$
the sensitivity to common and rare species. For $m <1$ the index reduces relative differences between abundant and rare species, while for $m>1$ exacerbates such differences, disproportionately favoring the most common species. Simpson's index is easily recovered for $m=2$, $H_2(P)=1- \sum_{i} P_i^2$ and Shannon entropy $H_{1}(P)=-\sum_i P_i \log P_i$,  the unique index that weighs all species exactly by their frequencies, for $m \rightarrow 1$. 
% For $P:=(P_i)_{i=1}^{W}$ an unknown distribution on a finite (or countably infinite) alphabet with cardinality $W$ is given by
%Like for Shannon and Simpson indexes, for $K$ a finite number of species, $H_m(P)$ reaches its maximum $(m-1)^{-1}(1- K^{1-m})$ at equiprobability, which reduces to $(m-1)^{-1}$ for $K \rightarrow \infty$. 
Few years later, in 1988, C. Tsallis introduces (\ref{tsallis}) in statistical physics  as a subadditive generalization of Shannon entropy, thus satisfying for $m >1$
$$
H_m(P_{A B})= H_m(P_A) + H_m(P_B) + (1-m)H_m(P_A) H_m(P_B),
$$
for $P_A$ and $P_B$ the relative abundances of two non overlapping independent classifications and $m$ a parameter measuring the degree of deviation from additivity. Since then this index, known as {\it Tsallis entropy}, plays a significant role in non-extensive generalizations of statistical mechanics  (Tsallis, 2009) and finds application in fields in which complex phenomena exhibit a power-law behaviour, reflecting a hierarchical or fractal structure. See e.g. Martins {et al.} 2009, Vila { et al.},  (2011), Zhang { et al.} (2010) for applications in machine learning, document classification, image processing and neural signals analysis. In community ecology $H_m(P)$ was rediscovered  by Keylock  in 2005 as a concave generalization of Simpson's and Shannon's measures addressing the self-similar nature of species abundances, as well as the significant amount of complex interactions between species and individuals in ecological systems. For a thorough analysis of entropy-based indices in ecology and their interpretation as diversity measures by a transformation in effective number of species see Jost (2006) and Mendes { et al.} (2008).

The typical problem in estimating diversity indices from a finite set of experimental data is that relative abundances are {a priori} unknown, and replacing them by sample relative frequencies, as in the maximum likelihood approach,
%by using the so-called {\it plug in} estimators,  
produces negatively biased estimators, especially in biological communities where a large number of species has relatively small abundances and many of the  rare species remain unobserved (cf. e.g. Chao and Shen, 2003).  In this perspective the Bayesian approach to diversity estimation has been already advocated as a more suitable solution.
% both for Shannon entropy like for the generalized Tsallis entropy. 
%for small datasets since no unbiased estimator is available and the maximum likelihood estimator based on frequency counts, the so-called {\it plug in} estimator, is severely biased when many of the species remain unobserved. 
Under the hypothesis  of finite and known number of species a first result for Shannon entropy estimation under symmetric Dirichlet priors is in Gill and Joanes (1979).  Independently in 1995 Wolpert and Wolf provide  posterior first and second moments under uniform prior on the finite dimensional simplex  devising a technique to obtain analogous results for more general priors.  Under the same hypothesis, in the setting of information theoretical analysis of neural responses,  Nemenman {et al.} (2002, 2004) show that symmetric Dirichlet priors impose a too narrow prior on Shannon entropy and suggest to use as an alternative a specific mixture of those priors. As for Tsallis diversity a recent and exhaustive analysis of maximum likelihood estimation compared to computationally intensive estimation methods is in Butturi-Gomes {et al.} (2014) 
%in which the efficiency of the maximum likeloohd estimator is compared tA 
while, to the best of our knowledge, the unique Bayesian proposal is under uniform prior on the finite dimensional simplex (Holtse {et al.}, 1998). 

Here we present a fully general Bayesian nonparametric solution, under a large class of priors, to the problem of estimating the general index $H_m(P)$ when  the relative abundances of the species in the population are unknown, the number of species is unknown and possibly countably infinite and the size of the sample available from the population is small and so is the number of different species observed.  Preliminary contributions along those lines, providing explicit posterior mean and variance under the two-parameter extension of the Poisson-Dirichlet distribution (Pitman and Yor, 1997), are in Cerquetti (2012) and Archer {et al.} (2013), respectively for Simpson's and Shannon's index. 
% class of discrete distributions generalizing the law of the ranked atoms of the Dirichlet process priors. 
Here  not only we are able to generalize those results to the whole family of Tsallis indices, but we derive posterior moments under a large class of priors introduced in Gnedin and Pitman (2006), which extends 
%introduce a remarkable class of infinite random discrete distributions in one to one correspondence with a family of (consistent) infinite exchangeable random partitions of the positive integers,  
the family of two-parameter Poisson-Dirichlet distributions while conserving the {Gibbs} product form of the corresponding exchangeable partition probability function. This class, which is commonly referred to as the {Gibbs priors} class, has become extremely popular in modern Bayesian nonparametrics as a tractable generalization of the Dirichlet process prior (Ferguson, 1973). Related Bayesian estimation in species sampling applications has been devised in Lijoi {et al.} (2007, 2008) and largely addressed in Favaro {et al.} (2009, 2012b, 2013). See also Cerquetti (2011, 2013, 2013b). Nevertheless, until now, the focus has been on posterior predictive species richness and discovery probability estimation. Here we provide the first results in this general setting for a large family of measures of diversity. 
% estimation complete the picture of Bayesian nonparametric estimation of species diversity by proposing posterior estimation for Patil-Tailllie 
Notice that despite (\ref{tsallis}) is defined for any $m > 0$, we will restrict to the case $m \in N$.
%which measures diversity by the number of different species in the population. ....... UFFF    
%Here we provide a full analysis of diversity as measured by  we complement and complete t face Bayesian estimation for 
We start by briefly introducing the Gnedin--Pitman class and some of its main properties. 
\subsection{Gnedin--Pitman priors}
Given an infinite random discrete distribution $P=(P_i)_{i \geq 1}$, 
%and $(X_n)_{n \geq 1}$ an infinite exchangeable sequence driven by $P$,
then the law of the infinite exchangeable random partition $\Pi_n=\{A_1, \dots, A_k\}$ of $[n]$, for $n \geq 1$, induced by sampling from $P$ 
%i.e. such that 
%By Kingman's correspondence, (\cite{kin78} Kingman, 1978), every $\Pi$ has the same distribution as one generated by an infinite exchangeable sequence $(X_n)$ driven by some random discrete probability measure $P=(P_i)_{i \geq 1}$: if 
%$i$ and $j$ belong to the same block iff $X_i=X_j$, 
is given by
\begin{equation}
\label{kingm}
p(n_1, \dots, n_k)=\sum_{(i_1, \dots, i_k)} {E} \left[  \prod_{j=1}^k P_{i_j}^{n_j} \right],
\end{equation}
where $n_i=|A_i|$,  $(i_1, \dots, i_k)$ ranges over all ordered $k$-tuples of distinct positive integers and $(P_i)_{i \geq 1}$ is any rearrangements of the ranked atoms $(P_i^{\downarrow})_{i \geq 1}$ of $P$. See Pitman (2003, 2006) for exhaustive accounts on exchangeable random partitions.
%t partitions arise by sampling 
%from {\it species sampling models} (cf. \cite{pit96}) 
%almost surely discrete probability measure $P=\sum_{i=1}^{\infty} P_i \delta _{X_i}(\cdot)$, by the exchangeable equivalence relation $i \approx j \iff X_i=X_j$.   )METTO QUI O METTO DOPO????) 
Gnedin--Pitman priors are the largest class of infinite random  discrete distributions with corresponding 
%exchangeable partition{\it Gibbs} partitions  (\citep{gnepit06}) 
%are the largest class of infinite 
exchangeable partition probability function (EPPF) (\ref{kingm}) in the {Gibbs} product form 
\begin{equation}
\label{EPPFgibbs}
p_{\alpha, V}(n_1, \dots, n_k)= V_{n,k} \prod_{j=1}^{k} (1 -\alpha)_{n_j-1},
\end{equation}
for $\alpha \in (-\infty, 1)$ and $V=(V_{n,k})$ weights satisfying the backward recursive relation 
\begin{equation}
\label{backR}
V_{n,k}=(n -k\alpha)V_{n+1, k} + V_{n+1, k+1},
\end{equation} 
where $V_{1,1}=1$ and $(x)_y=(x)(x+1)\cdots(x+y-1)$ is the usual notation for rising factorials. 
Given $(n_1, \dots, n_k)$ the multiplicities of the first $k$ species observed in a random sample of size $n=\sum_i n_i$,  the probabilities to observe the $j$th {old} species or a {new} yet unobserved species at step $n+1$ are easily obtained from (\ref{EPPFgibbs}) and correspond respectively to $p_{\alpha, V}(n_j^{+})= {V_{n+1,k}} (n_j -\alpha)/{V_{n,k}}$ and $p_{\alpha, V}(n^{k+1})=V_{n+1, k+1}/V_{n,k}$.
%\begin{equation}
%\label{onestep}
%p_{\alpha, V}(n_j^{+})= \frac{V_{n+1,k}}{V_{n,k}}(n_j -\alpha), \quad \mbox{} \quad
%p_{\alpha, V}(n^{k+1})=\frac{V_{n+1, k+1}}{V_{n,k}}.
%\end{equation}
By Theorem 12 in Gnedin and Pitman (2006) each random partition belonging to the class (\ref{EPPFgibbs}) is a probability mixture of extreme partitions, namely: finite symmetric Dirichlet partitions  for $\alpha < 0$, Ewens $(\theta)$ partitions (Ewens, 1972) for $\alpha=0$, and Poisson-Kingman conditional partitions driven by the stable subordinator (Pitman, 2003) for $\alpha \in (0,1)$. Therefore each element of the Gnedin-Pitman priors family can be obtained by mixing corresponding extreme random discrete distributions.
% namely finite Dirichlet distribution, over the number of species, Dirichlet process prior over the parameter $\theta$ and Poisson-Kingman priors over the total mass. 
%Besides Dirichlet priors and two-parameter Pitman Yor priors, popular Gnedin Pitman priors in Bayesian nonparametrics includes, two-parameter Gnedin-Fisher priors (Gnedin, 2010), normalized Stable priors,  and priors arising by exponential tilting of the Stable density, like normalized Generalized Gamma and normalized Inverse Gaussian priors (see Pitmna, 2003, Lijoi et al. 2005, 2007, Cerquetti, 2007). 
%The {\it two parameter $(\alpha, \theta)$ Poisson-Dirichlet} model (\cite{pit95}, \cite{pityor97}), which is well-known to arise as a polynomial mixture of Poisson............
Mathematical tractability characterizes some specific subfamilies of (\ref{EPPFgibbs}) that are typically the most widely implemented in the modern Bayesian nonparametric literature. Consequently first we provide closed form general expressions for prior and posterior moments of Tsallis measures under general Gnedin--Pitman priors and, as a by product, of Shannon's and Simpson's diversities. Then we derive explicit corresponding formulas  under {two-parameter $(\alpha, \theta)$ Poisson-Dirichlet priors}, and its particular cases, normalized $(\alpha)$ Stable and $(\theta)$ Dirichlet priors, and under the  { two-parameter Gnedin-Fisher priors} (Gnedin, 2010, Cerquetti, 2011b).  
%TOGLIAMO??? and  is given by
%$$
%p_{\gamma}(n_1, \dots, n_k)=\frac{(\gamma)_{n-k}(1)_{k-1}(1-\gamma)_{k-1}}{(1)_{n-1}(1+\gamma)_{n-1}}\prod_{j=1}^k (2)_{n_{j-1}}.
%$$
As for the class of {exponentially tilted} Poisson-Kingman priors driven by the stable subordinator, and in particular for {normalized Inverse Gaussian} and {normalized generalized Gamma} priors,  which belong to the Gnedin--Pitman class for $\alpha \in (0,1)$, (cf. Pitman, 2003), the explicit derivation of prior and posterior moments of Tsallis diversity and Shannon entropy  from our general results will be the topic of a future paper. 

We stress that the technique adopted here allows to obtain the full sequence of prior and posterior moments of $H_m(P)$  under $(\alpha, V)$ Gnedin--Pitman priors. Nevertheless we give explicitly just the first three moments for $H_m(P)$ and the first two moments for $H_1(P)$. 
%This is enough for Bayesian point estimation under quadratic loss function, by posterior mean, to bound the probability of deviation by Chebyschev's inequality, and to study the symmetry of the posterior distribution.
\section{Prior analysis}
\subsection{Tsallis diversity}
One of the main problem in Bayesian estimation of Shannon entropy under symmetric Dirichlet distributions on the relative abundances, is that this class induces an extremely concentrated distribution on the prior belief, with variance that vanishes as the number of species becomes large (cf. Nemenman { et al.} 2002, 2004). The very same problem arises under two parameter  $(\alpha, \theta)$ Poisson-Dirichlet priors on the infinite dimensional simplex for $\theta$ large (cf. Archer {et al.} 2013). 
%In both cases to obtain a prior which generates a nearly uniform distribution of Shannon entropy the authors consider a prior by mixing symmetric Dirichlet priors o.... In order to provide the tools 
To check for analogous pathological behaviours induced by priors belonging to the Gnedin-Pitman class on the distribution of $H_m(P)$ for $m \geq 2$, we first derive explicit closed form expressions for the first three {prior} moments. 
%We then specialize formulas to Simpson's index and Shannon entropy particular cases.

\begin{theorem} Given a random discrete distribution $P=(P_i)_{i \geq 1}$, distributed according to a $(\alpha, V)$ Gnedin-Pitman model with EPPF (\ref{EPPFgibbs}), then %first, second and third prior moments of $H_m(P)$ for $m \in N$ are given respectively by
%for $S_m=\sum_{j=1}^\infty P_j^m$, 
%prior mean and variance of Tsallis entropy $H_m(P)$ 
\begin{equation}
\label{mean_tsapr}
E_{\alpha, V}(H_m(P))= (m-1)^{-1}[1 -V_{m,1}(1-\alpha)_{m-1}], 
\end{equation}
\begin{eqnarray}
\label{sec_tsapr}
E_{\alpha, V}[(H_m(P))^2]&=&(m-1)^{-2} \left[1 +V_{2m, 1}(1-\alpha)_{2m-1}\right.\nonumber\\
&+&\left. V_{2m,2}((1-\alpha)_{m-1})^2 -2V_{m,1}(1-\alpha)_{m-1}\right]
\end{eqnarray}
%\begin{equation}
%\label{var_tsapr}
%var_{\alpha, V}(H_m)= (m-1)^{-1} 
%\left[V_{2m,1}(1-\alpha)_{2m-1}+ V_{2m,2}((1-\alpha)_{m-1})^2 - %(V_{m,1}(1-\alpha)_{m-1})^2\right].
%\end{equation}
and
\begin{eqnarray}
\label{mom_ter}
E_{\alpha, V}[(H_m(P))^3]&=&(m-1)^{-3} \left[1- 3V_{m,1}(1-\alpha)_{m-1} \right.\nonumber\\ 
&+&3[V_{2m,1}(1-\alpha)_{2m-1}+V_{2m,2}((1-\alpha)_{m-1})^2] \nonumber\\
&-&\left.V_{3m, 1}(1-\alpha)_{3m-1}+ 3V_{3m,2}[(1-\alpha)_{m-1}(1-\alpha)_{2m-1}]\right.\nonumber\\ &+&\left.V_{3m,3}((1-\alpha)_{m-1})^3 \right].
\end{eqnarray}
\end{theorem}
Specializing (\ref{mean_tsapr}),  (\ref{sec_tsapr}) and (\ref{mom_ter}) for $m=2$, and exploiting the backward recursion (\ref{backR}), corresponding prior moments for Simpson's index $H_2(P)=(1 -\sum_j P_j^2)$ easily follow. 
%\begin{proof}
%See the Appendix.
%\end{proof}
%{\bf Remark 1.} Level of uncertainty... $V_{m,1}(1-\alpha)_{m-1}$ is the probability of  COME STA IN RELAZIONE STA COSA CON LA PROBABILITA? DI OLD O DI NEW:::: ETC???? \\\\
\subsection{Shannon entropy}
As for Shannon entropy $H_1(P)$, the next theorem generalizes to the entire Gnedin-Pitman family the results on prior first and second moments under symmetric Dirichlet priors and under two-parameter Poisson-Dirichlet priors, already obtained respectively in Gill and Joanes (1979) and Archer {et al.} (2013).
\begin{theorem}\it Given a random discrete distribution $P=(P_i)_{i \geq 1}$, distributed according to a $(\alpha, V)$ Gnedin-Pitman model with EPPF (\ref{EPPFgibbs}), then 
%first prior moment of $H_1(P)$ is given by
\begin{eqnarray}
\label{SHApr_mean}
E_{\alpha, V}(H_1)&=& - \lim_{m \rightarrow 1} \frac{\partial}{\partial m}V_{m,1}-\psi_{0}(1-\alpha),
\end{eqnarray}
where $\psi_0$ is the {\it digamma } function. For $V^*_{r,s}=\lim_{m \rightarrow 1} \frac{\partial}{\partial m} V_{rm,s} $ and
$V^{**}_{r,s}=\lim_{m \rightarrow 1} \frac{\partial^2}{\partial m^2} V_{rm,s} $ and $\psi_1(\cdot)$ the {\it trigamma function} then 
%the prior second moment results
\begin{eqnarray}
\label{priSHAgibbs2}
E_{\alpha, V}[(H_1)^2]&=& \frac{1}{2} \left[4 \psi_0(2-\alpha)V^*_{2,1}+ 4 V_{2,1}\psi_0(2-\alpha)^2 + 4V_{2,1}\psi_1(2-\alpha)\right.\nonumber\\
&+&\left.V^{**}_{2,1} + 4 \psi_0(1-\alpha)V^*_{2,2}+4V_{2,2} \psi_0(1-\alpha)^2 +2 V_{2,2}\psi_1(1-\alpha) + V^{**}_{2,2}\right]\nonumber\\ 
&-&[2 \psi_0(1-\alpha)V^*_{1,1}+V_{1,1}\psi_0(1-\alpha)^2+V_{1,1}\psi_1(1-\alpha)+V^{**}_{1,1}].
\end{eqnarray}
\end{theorem}
The derivation of the third moment follows along the same lines. For brevity we omit here the explicit derivation.
\begin{remark} If the distributions of the size-biased atoms $(\tilde{P}_j)_{j \geq1}$
%the relative abundances in order of appearance in a size-biased sampling,
of the specific Gnedin--Pitman prior are known explicitly, like e.g. when a stick-breaking construction of the kind $\tilde{P}_1=V_1, \tilde{P}_j=V_j \prod_{i=1}^{j-1}(1-V_i)$ has been devised, then first and second prior moments of Shannon entropy can also be obtained through the same route adopted in Archer {et al.} (2013) under two-parameter Poisson-Dirichlet priors. This is also the case, for example, of {normalized Inverse Gaussian} priors, for which the stick breaking construction has been recently obtained in Favaro {et al.} (2012). Notice in fact that in the proof of Theorem 2
$$
\lim_{m\rightarrow 1}\frac{\partial^2}{\partial m^2}E[(S_m)^2]= \lim_{m \rightarrow 1}2E\left[\frac{\partial}{\partial m}S_m\right]^2 + \lim_{m\rightarrow 1}2E\left[\frac{\partial ^2}{\partial m^2}S_m\right].
$$
Now 
$$
\lim_{m \rightarrow 1}2E\left[\frac{\partial}{\partial m}S_m\right]^2 = \lim_{m \rightarrow 1}2E\left[\frac{\partial}{\partial m}\sum_j P_j^m\right]^2 = 
$$
$$
=\lim_{m \rightarrow 1} 2 E[\sum_j P_j^m \log P_j]^2=2 E\left[\sum_j (P_j \log P_j)^2 + 2 \sum_{i \neq j}P_i P_j \log P_i \log P_j\right],
$$
and, by properties of size-biased distributions (see e.g. Pitman, 1996, 2003), for $\tilde{P_1}$ and $\tilde{P}_2$ the first and second size-biased atoms,
$$
E[\sum_i (P_j)^2 (\log P_j)^2]=E[\tilde{P}_1 (\log \tilde{P}_1)^2]
$$
and 
$$
2  E\sum_{i \neq j} P_i P_j (\log P_i)(\log P_j)]= 2 E(\log \tilde{P}_1 \log \tilde{P}_2 (1 -\tilde{P}_1)).
$$
Additionally 
$$
\lim_{m \rightarrow 1}2E\left[\frac{\partial^2}{\partial m^2}(S_m)\right]= \lim_{m \rightarrow 1}2E\left[\frac{\partial^2}{\partial m^2}(\sum_j P_j^m)\right]=
$$
$$ =\lim_{m \rightarrow 1}2E[\sum_j P_j^m(\log P_j)^2]= 
2E[\sum_j P_j(\log P_j)^2]= 2E[(\log \tilde{P_1})^2].
$$
Nevertheless the implementation of this technique could be difficult when raw and mixed moments of $\log (\tilde{P}_1)$ and $\log (\tilde{P_2})$  are not known in closed form. The results in Theorem 1 and Theorem 2, while being valid for the entire Gnedin-Pitman class, 
%generalize the results in Archer et al. (2013) to the entire Gibbs priors family (\ref{}).
allow to derive prior moments for Shannon entropy directly by the Gibbs weights. 
\end{remark}
\subsection {Examples}
\begin{figure}
\centering
\resizebox{14.5cm}{!}{
\begin{tabular}{rrrr}
\includegraphics{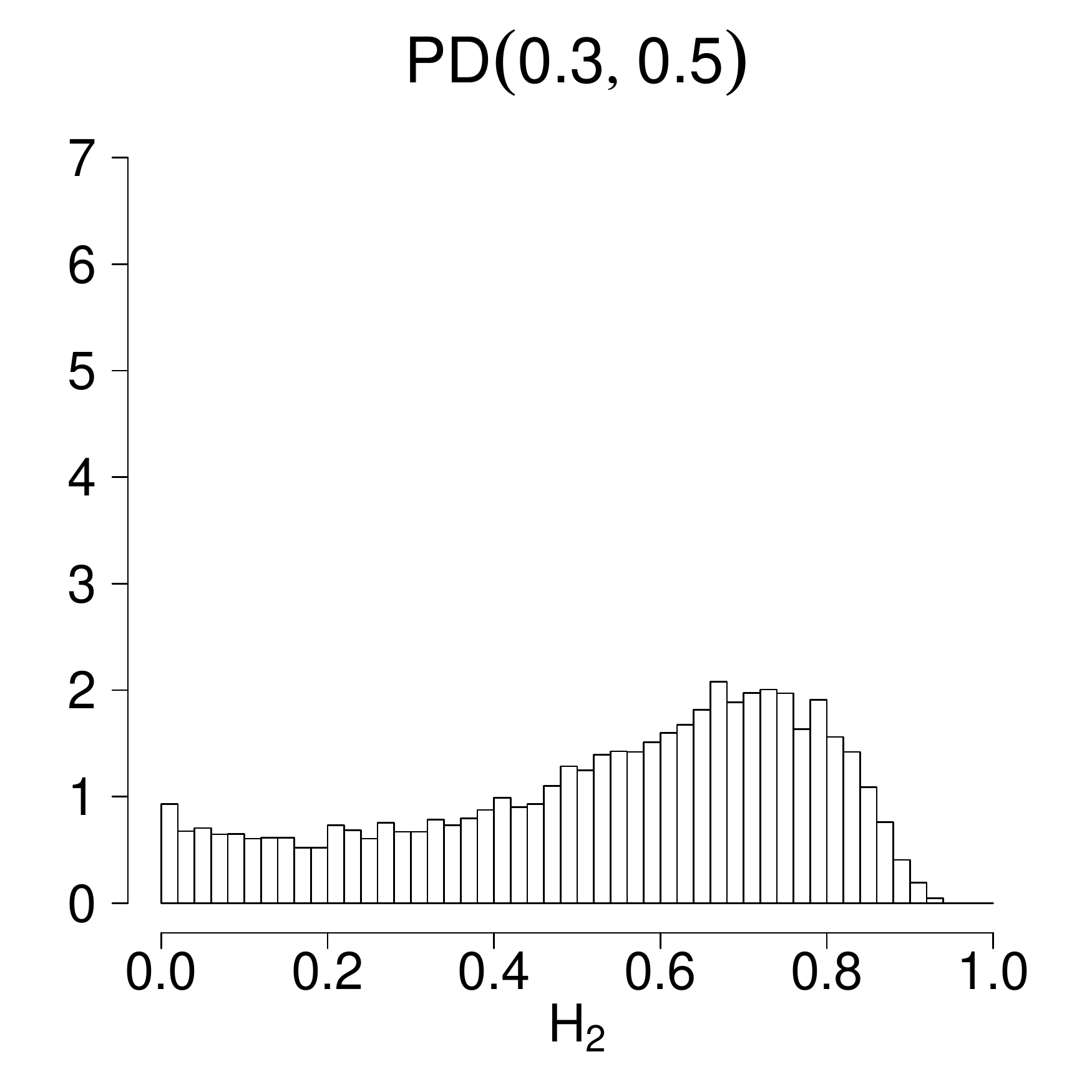}&
\includegraphics{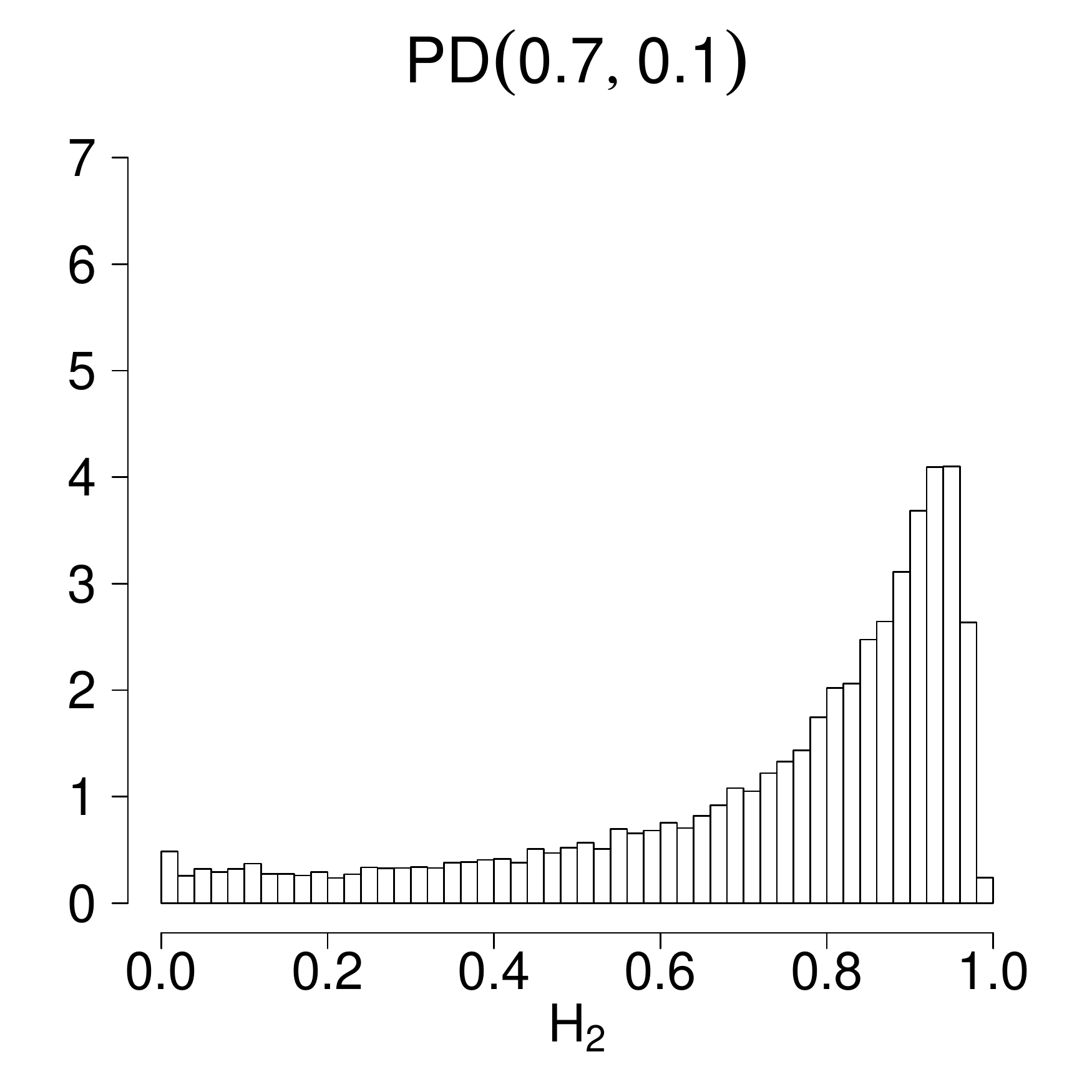}&
\includegraphics{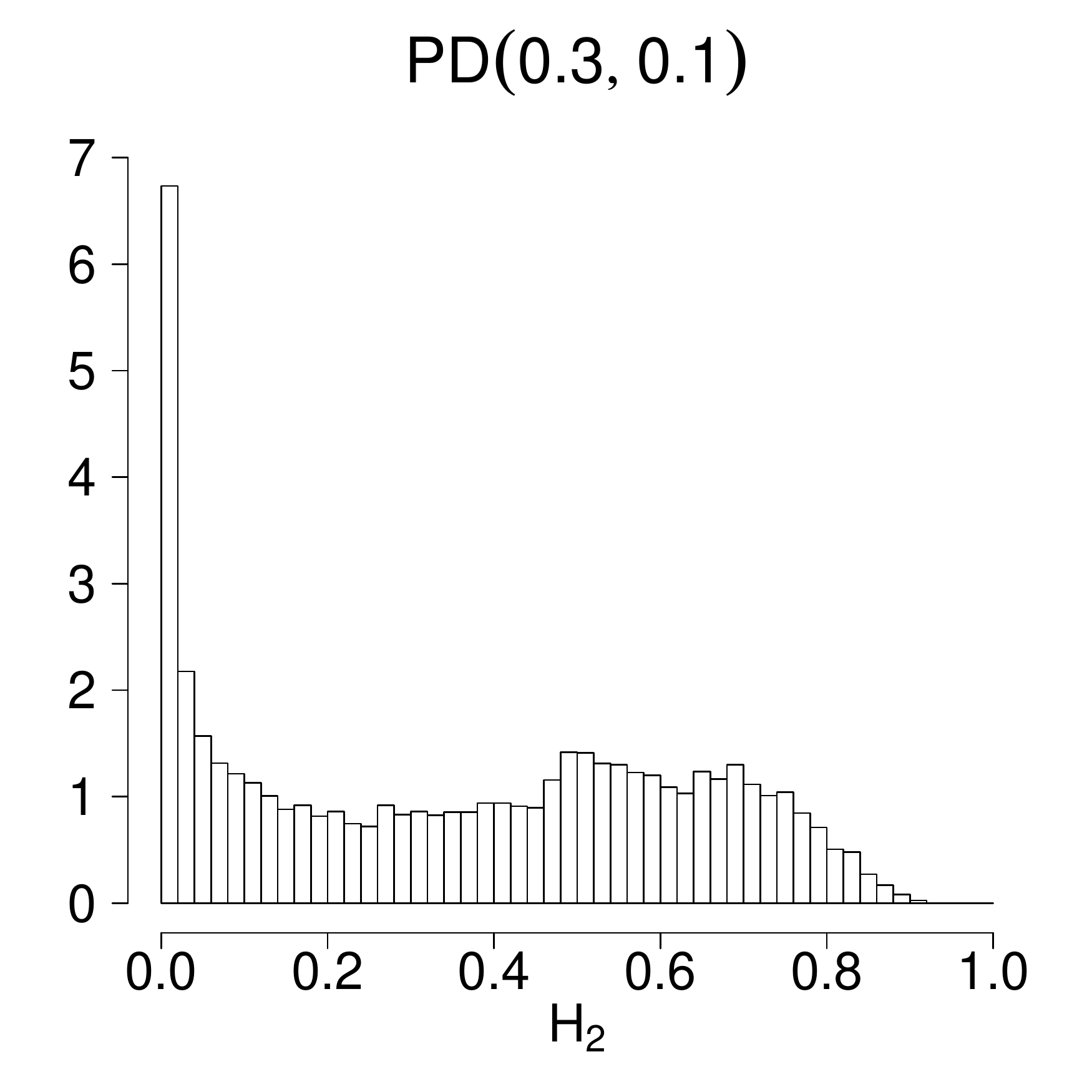}&
\includegraphics{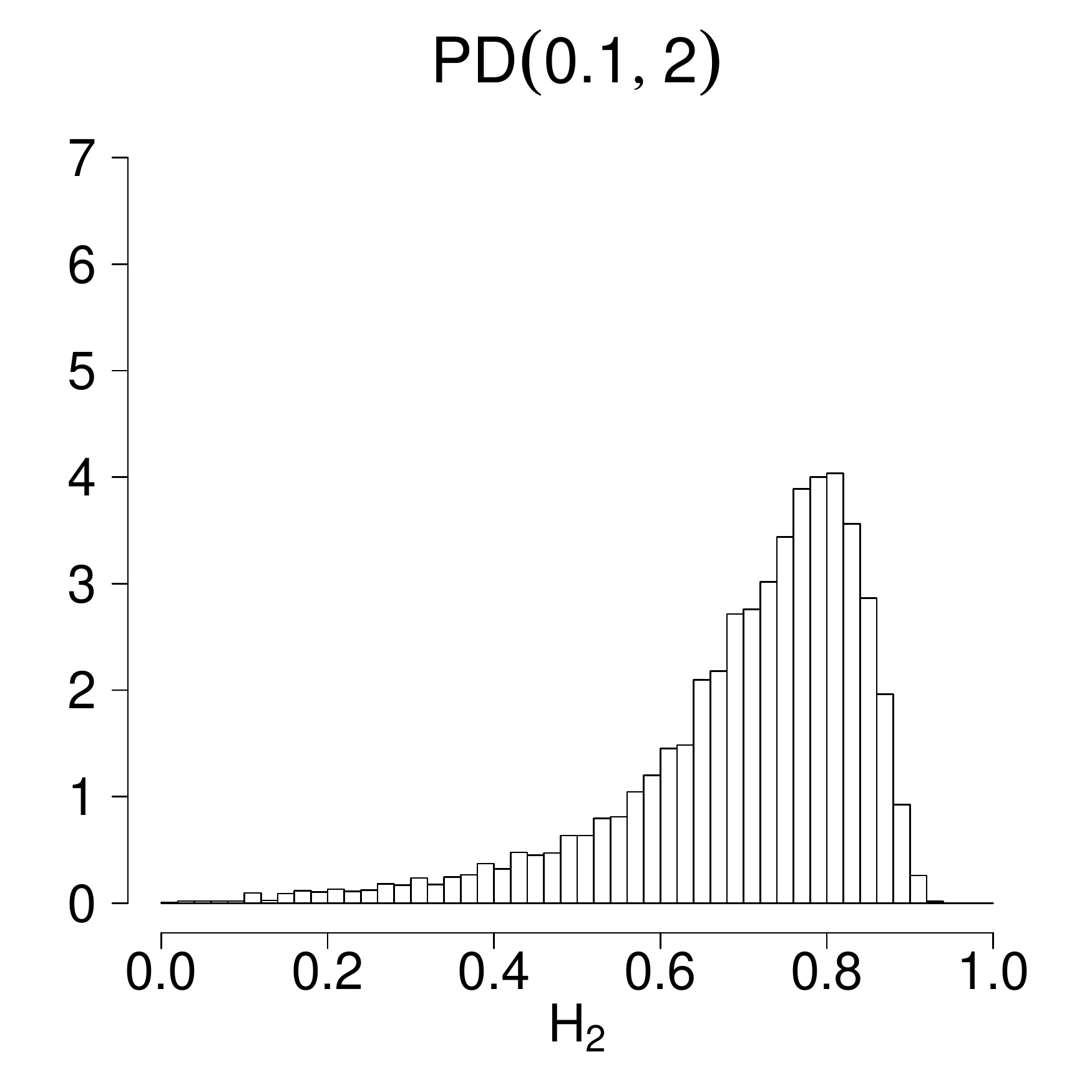}\\
\end{tabular}
}
\caption{Sampled prior distributions of $H_2$ under some $PD(\alpha, \theta)$ priors. $N=5000$, number of atoms truncated at $1000$.}
\label{histprior}
\end{figure}

\begin{example} [Two parameter Poisson-Dirichlet $(\alpha, \theta)$ priors] For $\alpha \in (0,1)$ and $\theta > -\alpha$ the {\it two parameter $(\alpha, \theta)$ Poisson-Dirichlet} model (Pitman, 1995, Pitman and Yor, 1997) has EPPF in the Gibbs form
\begin{equation}
\label{EPPFpd}
p_{\alpha, \theta}(n_1, \dots, n_k)= \frac{(\theta +\alpha)_{k-1 \uparrow \alpha}}{(\theta +1)_{n-1}} \prod_{j=1}^k (1 -\alpha)_{n_j-1},
\end{equation}
%\begin{equation}
%\label{2par}  
%p_{\alpha, \theta}(n_1, \dots, n_k)= 
%$V_{n,k}={(\theta +\alpha)_{k-1 \uparrow \alpha}}/{(\theta +1)_{n-1}}$,
% \prod_{j=1}^{k} (1 -\alpha)_{n_j-1}, 
%\end{equation}
for $(x)_{y \uparrow \alpha}=x(x+\alpha)\cdots(x+(y-1)\alpha)$ generalized rising factorials. %generalized rising factorials. and $(x)_y=(x)_{y\uparrow 1}$ . 
EPPFs of Dirichlet $(\theta)$  priors and normalized $(\alpha)$ Stable priors arise respectively for $\alpha=0$ and $\theta=0$ in (\ref{EPPFpd}).
An application of (\ref{mean_tsapr}) and (\ref{sec_tsapr}) yields
\begin{equation}
\label{priormeantsallisAT}
E_{\alpha, \theta}(H_m)= \frac{1}{m -1} \left(1 - \frac{(1-\alpha)_{m-1}}{(1+\theta)_{m-1}}\right)
\end{equation}
and 
\begin{equation}
\label{priorvariancetsallisAT}
 var_{\alpha, \theta} (H_m)= (m -1)^{-2} \left\{\frac{[(1-\alpha)_{m-1}]^2}{(1+\theta)_{m-1}} \left[\frac{(\theta +\alpha)}{(\theta+m)_m}- \frac{1}{(1+\theta)_{m -1}}\right] +\frac{(1-\alpha)_{2m -1}}{(1+\theta)_{2m-1}}\right\}.
\end{equation}
Simpson's index prior mean and variance follow by (\ref{priormeantsallisAT}) and (\ref{priorvariancetsallisAT}) for $m=2$ (see Cerquetti, 2012).
%$$
%{E}_{\alpha, \theta}(H_2)=\frac{\theta+\alpha}{1 +\theta}, \quad \mbox{and} \quad 
%var_{\alpha, \theta}(H_2)= \frac{(1 -\alpha)_3+ (\theta +\alpha)(1-\alpha)^2}{(\theta +1)_3} - \frac{(1-\alpha)^2}{(\theta +1)^2}.
%$$
Shannon entropy prior mean and variance, derived in Archer {et al.} (2013) exploiting the stick-breaking construction of the size-biased atoms, $\tilde{P}_1= V_1$ and $\tilde{P}_j= V_j \prod_{i=1}^{j-1}(1 -V_{i})$ for ${V}_j \sim Be(1-\alpha, \theta +j\alpha)$, follow by an easy application of (\ref{SHApr_mean}) and (\ref{priSHAgibbs2}): 
$$
E_{\alpha, \theta}(H_1)= \psi_0(\theta+1) - \psi_0(1-\alpha),
$$
%for $\psi_0=\frac{d}{dx} \log \Gamma(x)$ the {\it digamma function}, 
and
%As for the variance... (presa da Pillow, verificare con i miei conti che sia giusta e che non si possa ottenere piu' facilmente
$$
var_{\alpha, \theta} (H_1)= \frac{\theta +\alpha}{(\theta +1)^2(1- \alpha)}+ \frac{1-\alpha}{\theta +1} \psi_1(2-\alpha)-\psi_1(2+\theta).
$$
%for $\psi_1=\frac{d^2}{dx^2} \psi_0(x)$ is the {\it trigamma function}.
%\begin{example} [Dirichlet $(\theta)$ and normalized $(\alpha)$ Stable priors] 
Prior moments of Tsallis entropy, Simpson index and Shannon entropy under Dirichlet $(\theta)$ priors and normalized $(\alpha)$ Stable priors follow from the previous formulas respectively for $\alpha=0$ and $\theta=0$. 
\end{example}
Table~\ref{tabellePD}  illustrates numerically the prior choice on $H_m$ under two-parameter Poisson-Dirichlet priors for $m=2$ and $m=3$ and some combinations of $\alpha$ and $\theta$. Values are obtained by suitably applying formulas (\ref{priormeantsallisAT}) and (\ref{priorvariancetsallisAT}). For comparison purposes,  for $m=3$, the index the has been standardized by normalization with the maximum value $(m-1)^{-1}=1/2$. 
%First row  and column  correspond respectively to normalized Stable $(\alpha)$ prior and to Dirichlet $(\theta)$ prior cases. 
We refer the reader to Archer et al. (2013) for analogous prior analysis for Shannon entropy.  Notice that for $m=2$ under normalized Stable prior ($\theta=0$) the parameter $\alpha$ corresponds to the prior guess on the Simpson' index of diversity.  Increasing $\alpha$ for a given $\theta$ or increasing $\theta$ for a given $\alpha$ produces the same effect on the prior guess, which approaches the maximum value, and the same effect on the uncertainty, increasing the concentration of the prior around the mean.  This behaviour suggests,  as for Shannon index in Archer et al. (2013), that the choice of a two parameter Poisson-Dirichlet prior  for the generalized Tsallis diversity $H_m$ should be confined to combinations of small values of $\alpha$ and $\theta$.
%Under Dirichlet $(\theta)$ priors ($\alpha=0$), the prior guess on the diversity index increases with the value of  $\theta$.
%which is well-known to indicate the concentration of the prior guess mean value of the prior.
%\end{document}
\begin{table}[ht]
\def~{\hphantom{0}}
%\centerin
{\it Prior mean (left) and uncertainty (right, coefficient of variation) for standardized $H_m$ index under Poisson-Dirichlet priors for different values of $\alpha$ (columns) and $\theta$ (rows)}.
\resizebox{14.5cm}{!}{
\begin{tabular}{rrrrrrrrr}
&&&& {\large \it m=2} &&&\\\\
& $\alpha$ & 0 & 0.1 & 0.3 & 0.5 & 0.7 & 0.9 \\ 
 $\theta$ \\ %\hline
0 &&   & 0.100 & 0.300 & 0.500 & 0.700 & 0.900 \\ 
  0.1 && 0.091 & 0.182 & 0.364 & 0.545 & 0.727 & 0.909 \\ 
  0.5 && 0.333 & 0.400 & 0.533 & 0.667 & 0.800 & 0.933 \\ 
  1 && 0.500 & 0.550 & 0.650 & 0.750 & 0.850 & 0.950 \\ 
  1.5 && 0.600 & 0.640 & 0.720 & 0.800 & 0.880 & 0.960 \\ 
  2 && 0.667 & 0.700 & 0.767 & 0.833 & 0.900 & 0.967 \\ 
  4 && 0.800 & 0.820 & 0.860 & 0.900 & 0.940 & 0.980 \\ 
  10 && 0.909 & 0.918 & 0.936 & 0.955 & 0.973 & 0.991 \\ 
  12 && 0.923 & 0.931 & 0.946 & 0.962 & 0.977 & 0.992 \\\\  
\end{tabular}
\hspace{1cm}
\begin{tabular}{rrrrrrrrr}
&&&& {\large \it m=2} &&&\\\\
& $\alpha$ & 0 & 0.1 & 0.3 & 0.5 & 0.7 & 0.9 \\ 
 $\theta$ \\
0 &  && 1.732 & 0.882 & 0.577 & 0.378 & 0.192 \\ 
  0.1 && 1.753 & 1.176 & 0.733 & 0.506 & 0.339 & 0.175 \\ 
  0.5 && 0.676 & 0.586 & 0.447 & 0.338 & 0.239 & 0.128 \\ 
  1 && 0.408 & 0.369 & 0.300 & 0.236 & 0.171 & 0.094 \\ 
  1.5 && 0.291 & 0.267 & 0.222 & 0.178 & 0.132 & 0.073 \\ 
  2 && 0.224 & 0.207 & 0.174 & 0.141 & 0.105 & 0.059 \\ 
  4 && 0.109 & 0.102 & 0.088 & 0.073 & 0.055 & 0.031 \\ 
  10 && 0.036 & 0.034 & 0.030 & 0.025 & 0.019 & 0.011 \\ 
  12 && 0.028 & 0.027 & 0.023 & 0.020 & 0.015 & 0.009 \\\\
  % \hline
\end{tabular}
}
%\vspace{0.5cm}
\resizebox{14.5cm}{!}{ 
\begin{tabular}{rrrrrrrr}
&&&& {\large \it m=3} &&&\\\\
 & $\alpha$ & 0 & 0.1 & 0.3 & 0.5 & 0.7 & 0.9 \\ 
 $\theta$\\
0 &&  & 0.145 & 0.405 & 0.625 & 0.805 & 0.945 \\ 
  0.1 && 0.134 & 0.260 & 0.485 & 0.675 & 0.831 & 0.952 \\ 
  0.5 && 0.467 & 0.544 & 0.683 & 0.800 & 0.896 & 0.971 \\ 
  1 && 0.667 & 0.715 & 0.802 & 0.875 & 0.935 & 0.982 \\ 
  1.5 && 0.771 & 0.805 & 0.864 & 0.914 & 0.955 & 0.987 \\ 
  2 && 0.833 & 0.858 & 0.901 & 0.938 & 0.968 & 0.991 \\ 
  4 && 0.933 & 0.943 & 0.960 & 0.975 & 0.987 & 0.996 \\ 
  10 && 0.985 & 0.987 & 0.991 & 0.994 & 0.997 & 0.999 \\ 
  12 && 0.989 & 0.991 & 0.993 & 0.996 & 0.998 & 0.999 \\\\  
\end{tabular}
\hspace{1cm}
\begin{tabular}{rrrrrrrrr}
&&&& {\large \it m=3} &&&\\\\
& $\alpha$ & 0 & 0.1 & 0.3 & 0.5 & 0.7 & 0.9 \\ 
$\theta$\\
0 &&  & 1.705 & 0.838 & 0.525 & 0.326 & 0.155 \\ 
  0.1 && 1.739 & 1.147 & 0.689 & 0.454 & 0.288 & 0.139 \\ 
  0.5 && 0.647 & 0.549 & 0.401 & 0.287 & 0.190 & 0.093 \\ 
  1 && 0.371 & 0.328 & 0.252 & 0.186 & 0.126 & 0.063 \\ 
  1.5 && 0.250 & 0.224 & 0.176 & 0.131 & 0.089 & 0.045 \\ 
  2 && 0.182 & 0.164 & 0.130 & 0.098 & 0.067 & 0.033 \\ 
  4 && 0.072 & 0.065 & 0.052 & 0.040 & 0.027 & 0.014 \\ 
  10 && 0.015 & 0.013 & 0.011 & 0.008 & 0.006 & 0.003 \\ 
  12 && 0.010 & 0.009 & 0.008 & 0.006 & 0.004 & 0.002 \\\\
   %\hline
\end{tabular}
\label{tabellePD}
}
%\end{center}
\end{table}
Figure~\ref{histprior} shows four simulated prior distributions of Simpson's index under different combinations of $\alpha$ and $\theta$. Sampled values of $H_2(P)$ are obtained by resorting to the stick-breaking construction of the size-biased atoms of the $PD(\alpha, \theta)$ laws.  Notice that as for concentration, also the symmetry of the distribution dramatically changes with the values of the parameters and large values of $\alpha$ or $\theta$ produce increasing left skewed distributions.
\begin{figure}
\centering
\resizebox{14.5cm}{!}{
\begin{tabular}{rrrr}
\includegraphics{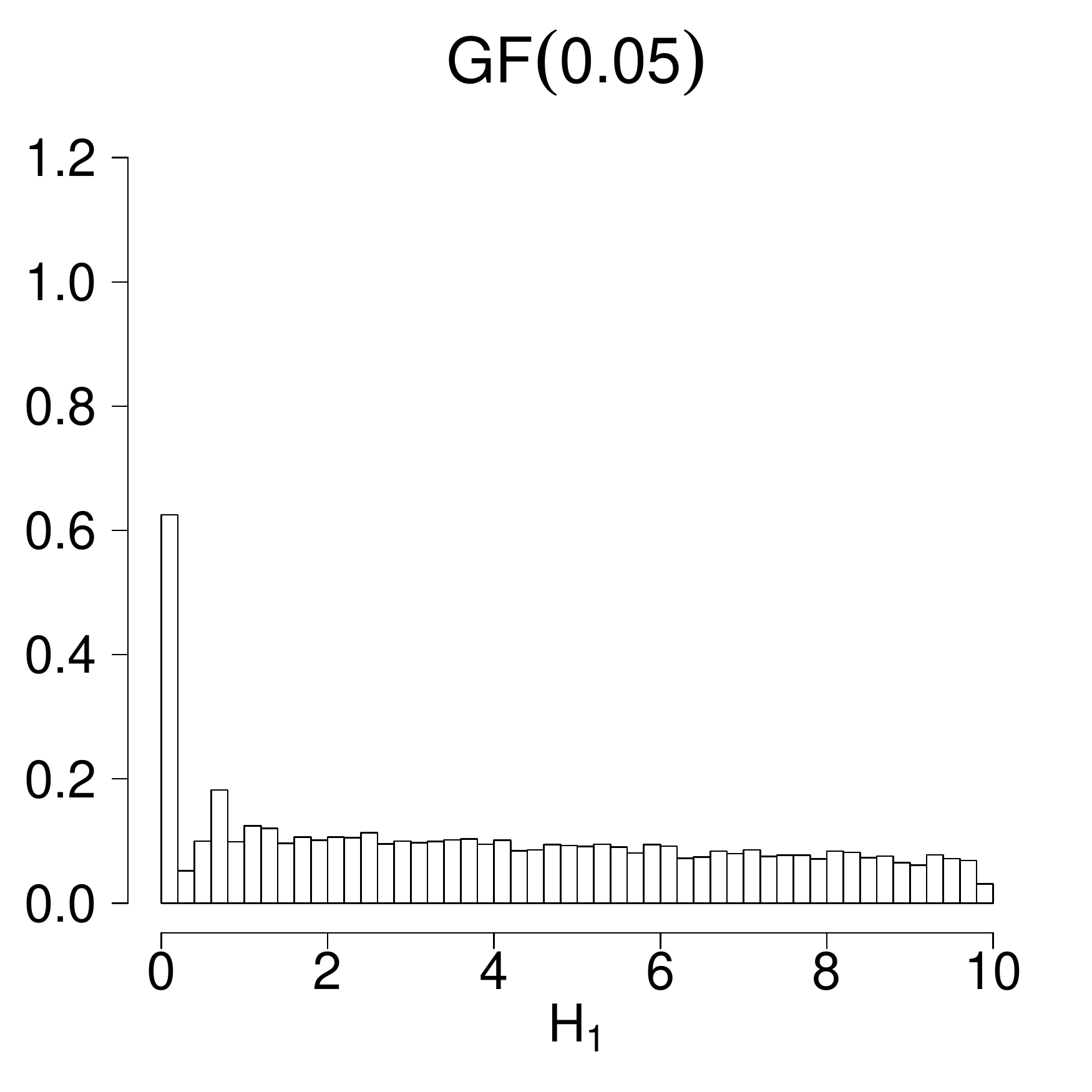}&
\includegraphics{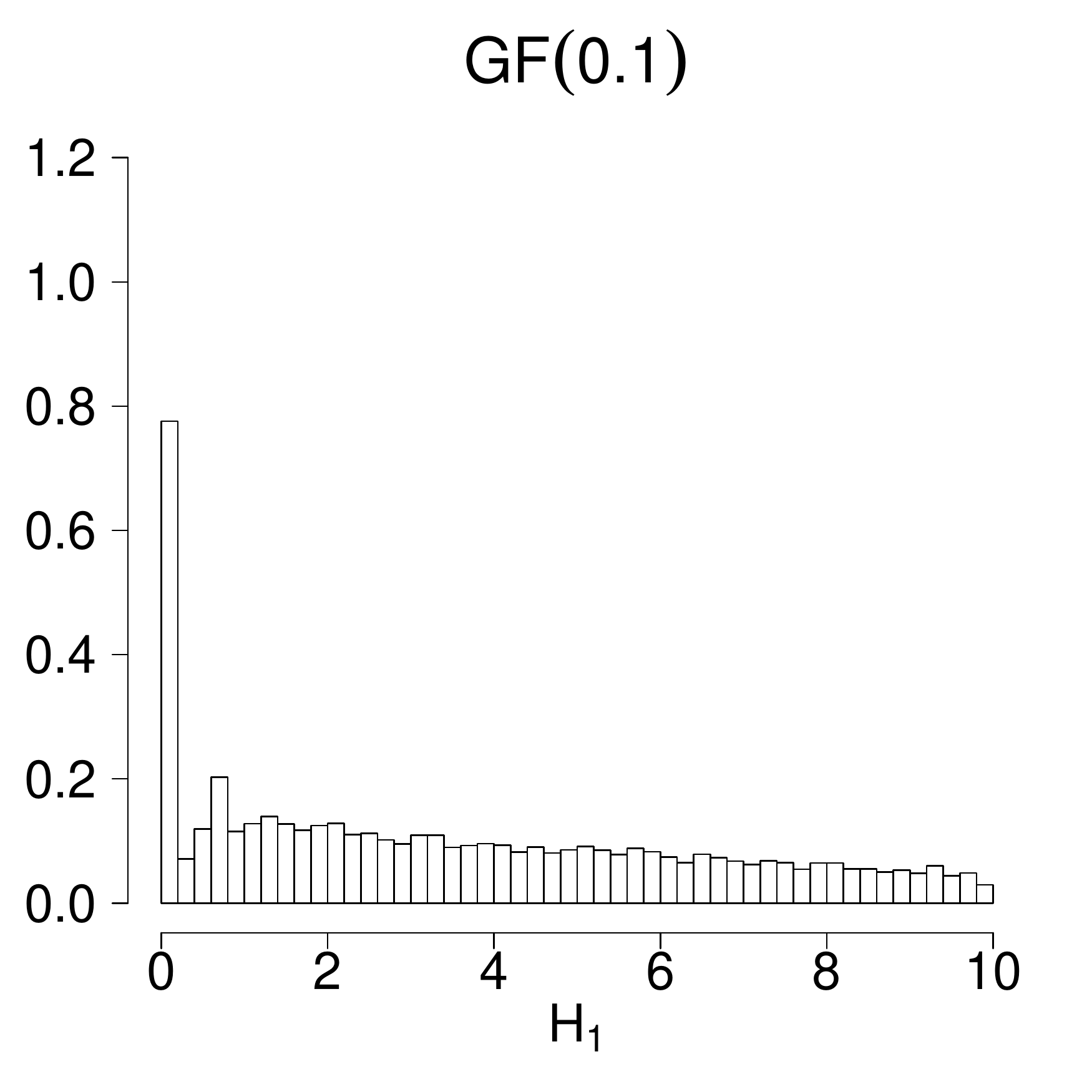}&
\includegraphics{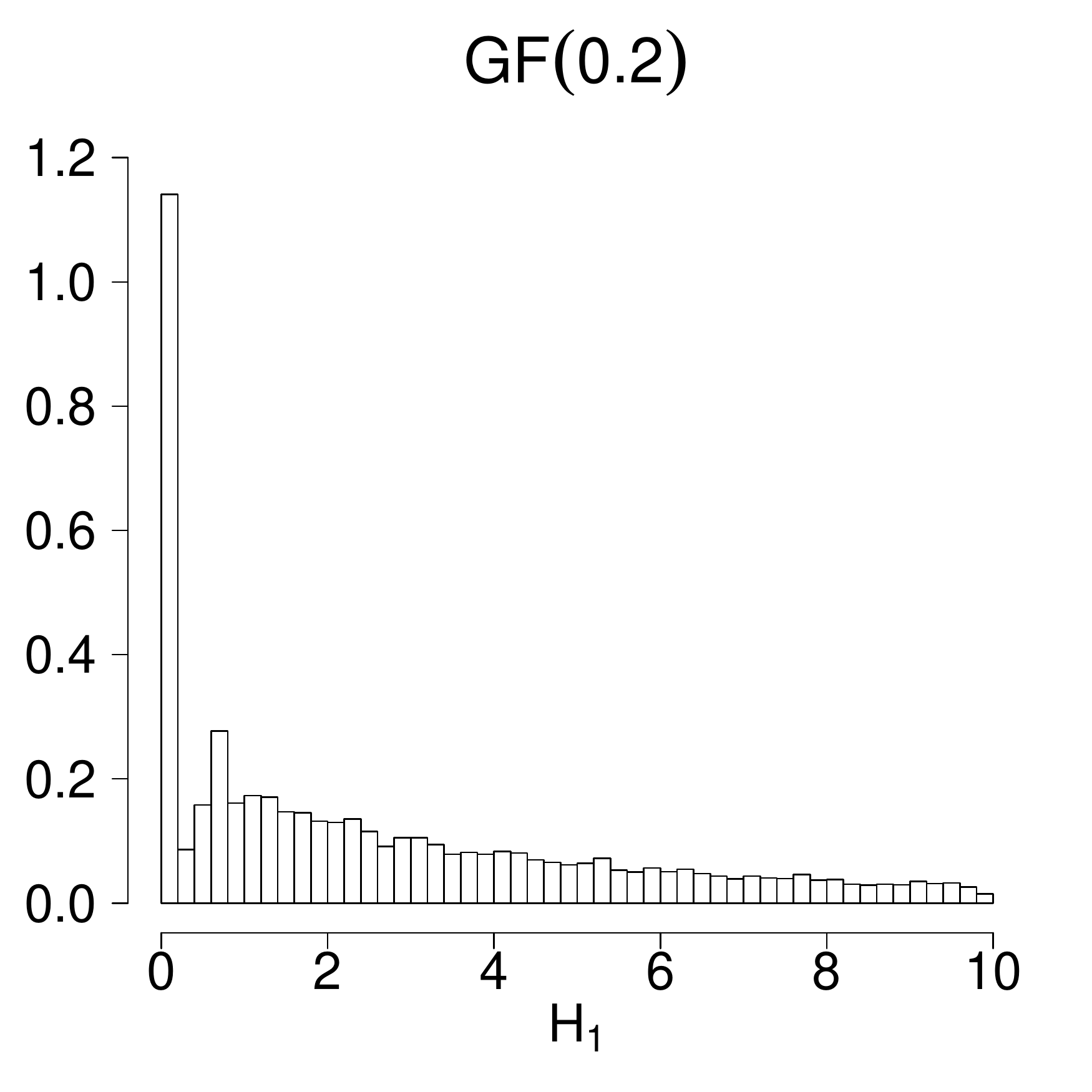}&
\includegraphics{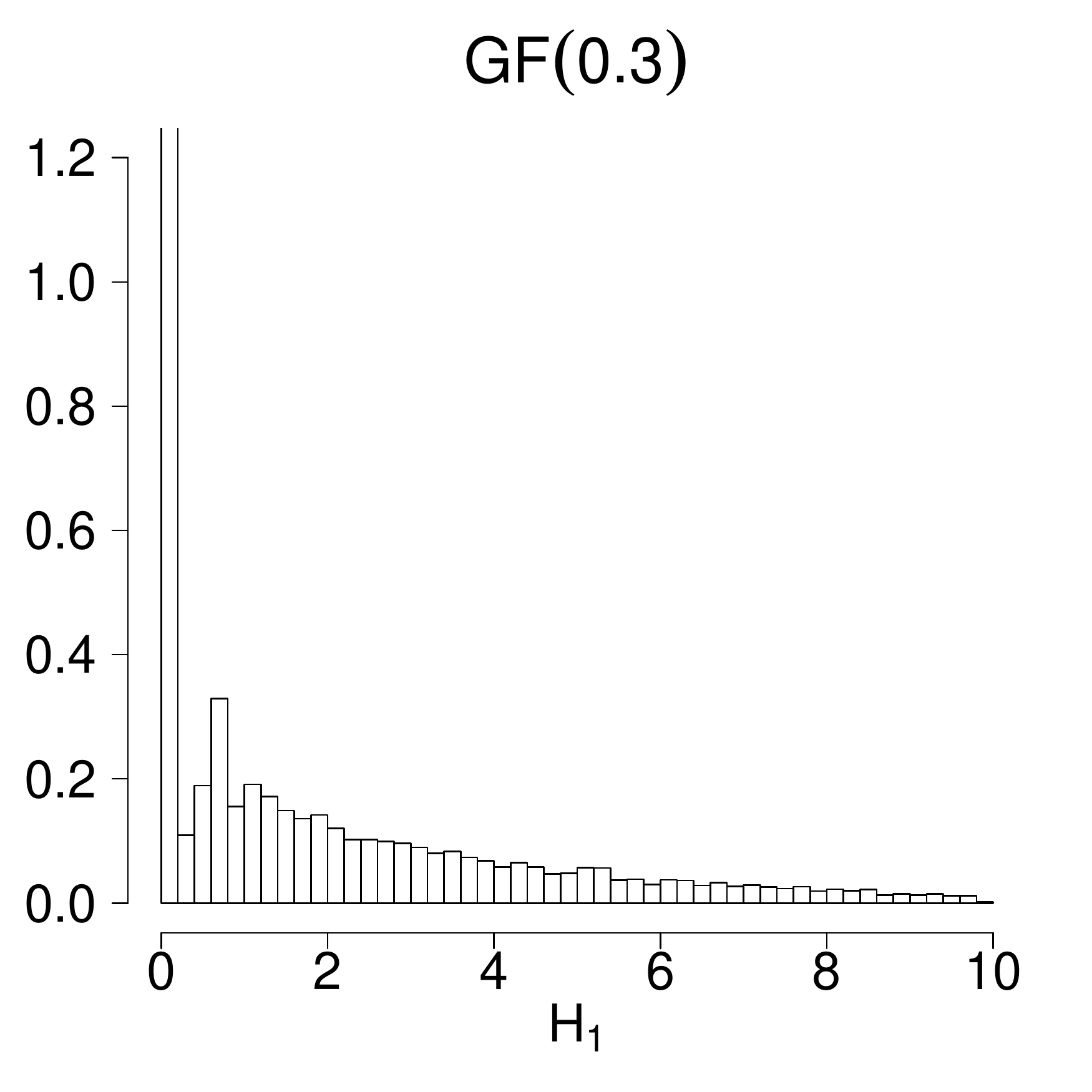}\\
\end{tabular}
}
\caption{Sampled prior distributions of Shannon entropy $H_1$ under one parameter Gnedin--Fisher priors for four different values of $\gamma$, $N=10.000$.}
\label{histpriorGNEDINFISHER}
\end{figure}
%Figure 1 shows simulated prior distributions of Simpson's index under four different combinations of $\alpha$ and $\theta$.  
%For $\theta >2$ the prior guess on $H_m$ is of a near uniform distribution of the population among species, and increases with the value of $\alpha$. 
 
%For high values of $\alpha$ and $\theta \geq 2$ the prior guess on $H_2$ tends to one, thus implying that the prior belief on the unknown abundances is concentrated on the hypothesis of uniform distribution of the population among species.
\begin{example} [Two parameter Gnedin-Fisher $(\psi, \gamma)$ priors] In 2010  Gnedin introduced another tractable two-parameter family of laws for random discrete distributions  belonging to the Gibbs class for $\alpha <0$. The model is obtained by mixing the uniform prior on the finite dimensional simplex with {\it shifted generalized Waring distributions} (Xekalaki, 1983) on the number of species. Here we adopt the parametrization devised in Cerquetti (2011) for its  specific tractability in the Bayesian nonparametric setting. The general form of the EPPF, for $\psi \in [0,1)$ and $0 < \gamma <\psi+1$,  is given by
\begin{equation}
\label{GF-EPPF2}
p_{\gamma, \psi}(n_1, \dots, n_k)=\frac{(\gamma)_{n-k}(1- \psi)_{k-1}(1-\gamma+\psi)_{k-1}}{(1+\psi)_{n-1}(1+\gamma- \psi)_{n-1}}\prod_{j=1}^k n_j!\nonumber
\end{equation}
and reduces to the one-parameter case for $\psi=0$.
An application of  (\ref{mean_tsapr}) and (\ref{sec_tsapr}) yields
\begin{equation}
\label{priormeanGF}
E_{\psi, \gamma}(H_m)=(m-1)^{-1}\left[1 -\frac{(\gamma)_{m-1}(2)_{m-1}}{(1+\psi)_{m-1} (1+\gamma -\psi)_{m-1}}\right],
\end{equation} 
with prior variance
\begin{eqnarray}
\label{vargnedprio}
var_{\psi, \gamma}(H_m)&=& (m-1)^{-2}\left[\frac{(\gamma)_{2m -1}(2)_{2m-1}+ (\gamma)_{2m-2}(1-\gamma +\psi)(1-\psi)[(2)_{m-1}]^2}{(1 +\psi)_{2m-1}(1+\gamma- \psi)_{2m-1}}\right. \nonumber\\
&-&\left.\frac{[(\gamma)_{m-1}(2)_{m-1}]^2}{[(1+\psi)_{m-1}(1+\gamma-\psi)_{m-1}]^2}\right].
\end{eqnarray}
Specializing (\ref{priormeanGF}) and (\ref{vargnedprio}) for $m=2$ corresponding formulas for Simpson's index easily follow.
As for Shannon entropy prior mean an application of (\ref{SHApr_mean}) yields 
%suitably mixing $PD(-1, -\xi)$ estimates. It follows that 
$$
E_{\gamma, \psi}(H_1)= -\psi_0(2)+\psi_0(1+\gamma-\psi)-\psi_0(\gamma)+\psi_0(1+\psi),
$$
which reduces to $E_{\gamma}(H_1)=  \gamma^{-1}{(1-\gamma)}
$ for the one-parameter case $(\psi=0)$.
%$$
%E_{\gamma, \psi}(H_E)= \gamma[\psi_0(\gamma+1)+\psi_0(2)-\psi_0(\psi+1)-\psi_0(1+\gamma-\psi)]
%$$
%and in fact corresponds to $\lim_{m \rightarrow 1} \frac{1}{m-1}(1 -\frac{\gamma m}{\gamma +m-1})$.
Prior second moment arises by calculating the limit for $m \rightarrow 1$ of the first and second partial derivatives of 
$$
V_{\xi m,i}^{\gamma, \psi}= \frac{(\gamma)_{\xi m-i}(1-\psi)_{i-1}(1-\gamma+\psi)_{i-1}}{(1+\psi)_{\xi m-i}(1+\gamma -\psi)_{\xi m-i}},
$$
with respect to $m$ and then applying (\ref{priSHAgibbs2}). For the sake of brevity we omit here the explicit formulas. 
Figure~\ref{histpriorGNEDINFISHER} shows some simulated prior distributions of Shannon entropy under Gnedin--Fisher one parameter prior for different values of $\gamma$. A peak at zero, which increases with $\gamma$, characterizes those priors since the Waring mixing distribution puts positive mass $\gamma$ at $\xi=1$,  the case of a population with a unique species. Nevertheless the priors on $H_1$ appear to be sufficiently flat and uninformative in  $(0, 10)$ for small values of $\gamma$, thus suggesting that those models may provide an interesting choice as priors for Shannon entropy, not suffering of the problem of concentration of the two parameter Poisson-Dirichlet model and its particular cases.  
\end{example}
\section{Bayesian nonparametric estimation}
\subsection{Tsallis diversity}
In the next theorem, which generalizes the results in Holste et al. (1998, Sect. 4),  we obtain the first three posterior moments of $H_m(P)$ thus providing what needed for Bayesian point estimation under quadratic loss function, for approximate interval estimation by Chebyschev's inequality and for studying the symmetry of the posterior distribution. Given the vector $(n_1, \dots, n_k)$ of the multiplicities of the first $k$ different species observed in order of appearance in a sample of size $n$, similarly to Theorem 1,  the results arise deriving the posterior moments of $S_m=\sum_{j=1}^\infty P_j^m$ under Gnedin-Pitman priors and then applying the binomial theorem. The details of the proof are in the Appendix.
%See Proposition A2 in the Appendix. 
%(\ref{primikappa}), (\ref{misti}) and (\ref{resto}) are enough to derive posterior moments of any order for $H_m(P)$ under Gnedin-Pitman priors in closed form. 
\begin{theorem} 
%Let $(X_i)_{i \geq 1}$ be an exchangeable random sequence of species observations driven by  
Let ${\bf n}=(n_1, \dots, n_k)$ be the multiplicities of the first $k$ species observed in a sample of size $n$, then, under a general $(\alpha, V)$ Gnedin-Pitman prior on the unknown relative abundances $(P_i)_{i \geq 1}$
%posterior first, second and third moments of the Patil-Taillie-Tsallis diversity $H_m(P)$ are given respectively by
\begin{eqnarray}
\label{post_tsagibbs}
E_{V, \alpha, \bf {n}}(H_m)= \frac{1}{m-1}\left[1 -\frac{V_{n+m,k}}{V_{n,k}}\sum_{j=1}^k (n_j -\alpha)_m - \frac{V_{n+m, k+1}}{V_{n,k}}(1-\alpha)_{m-1}\right],
\end{eqnarray}
%with posterior second moment
%\begin{equation}
%\label{post_tsavar}
%
\allowdisplaybreaks
%\begin{align*}
\begin{eqnarray}
\label{posttsa2}
E_{V, \alpha, {\bf n}}[(H_m)^2]&=& \left(\frac{1}{m-1}\right)^2 \left[1 +\frac{V_{n+2m, k}}{V_{n,k}}\left(\sum_j (n_j -\alpha)_{2m}+2\sum_{i \neq j} (n_j-\alpha)_m(n_i-\alpha)_m \right) \right.\nonumber\\
&+&\frac{V_{n+2m, k+2}}{V_{n,k}} [(1-\alpha)_{m-1}]^2 \nonumber\\
&+& \frac{V_{n+2m, k+1}}{V_{n,k}}[(1-\alpha)_{2m-1}+ 2 (1-\alpha)_{m-1}\sum_j (n_j-\alpha)_m ] \nonumber\\
&-&\left.2 \left(\frac{V_{n+m,k}}{V_{n,k}}\sum_{j=1}^k (n_j -\alpha)_m - \frac{V_{n+m, k+1}}{V_{n,k}}(1-\alpha)_{m-1}\right)\right]
\end{eqnarray}
%\end{align*}
and 
\allowdisplaybreaks
\begin{eqnarray}
\label{posttsa3}
E_{V, \alpha, {\bf n}}[(H_m)^3]=(m-1)^3 \{1 -3 E_{\alpha, V, {\bf n}}(S_m)+3E_{\alpha, V, {\bf n}}[(S_m)^2] - E_{\alpha, V, {\bf n}}[(S_m)^3]\} \nonumber
\end{eqnarray}
for
\allowdisplaybreaks
\begin{eqnarray}
%\begin{align*}
\label{postSm3}
E_{\alpha, V, {\bf n}}[(S_m)^3]&=&\frac{V_{n+3m, k}}{V_{n,k}}\left[\sum_{j=1}^k (n_j-\alpha)_{3m}+ 3 \sum_{i \neq j} (n_j -\alpha)_{2m}(n_i-\alpha)_{m}\right. \nonumber\\
&+& \left. 3! \sum_{i \neq j \neq h} (n_i -\alpha)_m (n_j -\alpha)_m (n_h-\alpha)_m\right]\nonumber\\
&+&   \frac{6V_{n +3m, k+1}}{V_{n,k+1}}(1-\alpha)_{m-1}[\sum_{j=1}^k (n_j -\alpha)_{2m}+ 2 \sum_{i \neq j} (n_j -\alpha )_m (n_i -\alpha )_m ]\nonumber\\
&+& \frac{3 V_{n+3m, k+2}}{V_{n,k}}[(1-\alpha)_{m-1}]^2 \sum_{j=1}^k (n_j -\alpha)_m \nonumber\\
&+& \frac{V_{n +3m, k+1}}{V_{n,k}}(1-\alpha)_{3m-1} + 3 \frac{V_{n+3m, k+2}}{V_{n,k}} (1-\alpha)_{m-1} (1-\alpha)_{2m-1}\nonumber\\
&+& \frac{V_{n+3m, k+3}}{V_{n,k}}[(1-\alpha)_{m-1}]^3.
\end{eqnarray}
%\end{align*}
\end{theorem}
Posterior moments of Simpson's index ($m=2$), generalizing the results in Cerquetti (2012) to the entire Gnedin-Pitman class, arise specializing (\ref{post_tsagibbs}), (\ref{posttsa2}) and (\ref{postSm3}) for $m=2$.
%\begin{proof}The results follow easily for $\xi=1$, $\xi=2$ and $\xi=3$ in (\ref{primikappa}), (\ref{misti}) and (\ref{resto}). 
%\end{proof}
%\begin{remark}
%The results in (\ref{post_tsagibbs}) and (\ref{posttsa2}) generalizes the results in Holtse et al. (1998), derived under uniform prior on the finite dimensional simplex,  to the entire Gnedin--Pitman family. In particular  equation (12) can be obtained from equation (\ref{postmean}) for $\alpha=-1$ and $\theta=M$, where $\theta$ is the number of different species.
%\end{remark}
\subsection{Shannon entropy}
%%%%%%%%%%%%%%%%%%%%%%%
%%PROPOSITION 2
%%%%%%%%%%%%%%%%%%%%%%%%%%%%
Similarly to Theorem 2, posterior estimation of Shannon entropy under general Gnedin--Pitman priors, which provides a substantial generalization of the results in Wolpert and Wolf (1995), Nemenmann {et al.} (2002, 2004) and Archer {et al.} (2013), is obtained by repeated application of H\^opital's rule to the results in Theorem 3. The following Proposition provides posterior mean and posterior second moment under a general $(\alpha, V)$ Gnedin--Pitman model. Notice that, despite it is possible to derive the closed form, the posterior second moment expression is extremely complex and not easy to be calculated in real data applications. Therefore, in Section 4, we illustrate the posterior uncertainty of $H_1(P)$ deriving the highest posterior density intervals (HPD) by the simulated posterior distributions. 
\begin{proposition}
Let ${\bf n}=(n_1, \dots, n_k)$ be the multiplicities of the first $k$ species observed in a sample of size $n$, then, under a general $(\alpha, V)$ Gnedin--Pitman prior on the unknown relative abundances,
\allowdisplaybreaks
\begin{align}
\label{postSHANgibbs}
E_{V, \alpha, {\bf n}}(H_1)&=
- \frac{\sum_j (n_j -\alpha)}{V_{n,k}}[\lim_{m \rightarrow 1}\frac{\partial}{\partial m}V_{m+n,k} +\psi_0(n_j-\alpha+1)V_{n+1,k}] \nonumber\\
&-\frac{1}{V_{n,k}}[\lim_{m \rightarrow 1} \frac{\partial}{\partial m} V_{n+m, k+1}+ \psi_0 (1-\alpha)V_{n+1, k+1}],
\end{align}
{\small {\normalsize and}
\allowdisplaybreaks
\begin{eqnarray}
\label{secpostsha}
 E_{\alpha, V, {\bf n}}[(H_1)^2]&=& \frac{1}{V_{n,k}} \left\{V_{n+2,k}^{**} [\sum_{j} (n_j -\alpha)_2 + 2 \sum_{i \neq j} (n_i-\alpha)(n_j -\alpha)]\right.\nonumber\\
&+&\left.V_{n+2,k}^{*}\{4 \sum_j (n_j -\alpha)_2 \psi_0(n_j -\alpha +2) \right.\nonumber\\
&+& \left.  2 \sum_{i\ \neq j} (n_j-\alpha)(n_i -\alpha)[\psi_0 (n_i -\alpha +1) + \psi_0 (n_j-\alpha+1 )]\}\right.\nonumber\\
&+& V_{n+2,k}\left\{ 4 \sum_{j}(n_j -\alpha)_2[\psi_0(n_j -\alpha+2)^2 +\psi_1(n_j -\alpha+2)^2]\right. \nonumber\\
&+&2 \sum_{i \neq j}(n_j -\alpha)(n_i -\alpha)\left[\psi_0(n_i -\alpha +1)\psi_0(n_j -\alpha +1)\right.\nonumber\\
&+&\left.\left.\psi_0(n_i -\alpha+1)^2 +\psi_0(n_j -\alpha +1)^2 +\psi_1(n_i-\alpha+1) +\psi_1(n_j- \alpha+1)\right]\right\}\nonumber\\
&+&4V_{n+2, k+2}^{**} + 8V_{n+2,k+2}^*\psi_0(1-\alpha) +2V_{n+2,k+2}(2\psi_0(1-\alpha)^2 +\psi_1(1-\alpha))\nonumber\\
&+&4 V_{n+2, k+1}^{**}[2 \sum_j (n_j-\alpha)+(1-\alpha)]\nonumber\\
&+&8 V_{n+2,k+1}^*\{\sum_j (n_j-\alpha)[\psi_0(n_j -\alpha+1)+\psi_0(1-\alpha)]+ (1-\alpha)\psi_0(2-\alpha)\}\nonumber\\
&+&2 V_{n+2, k+1}\left[\sum_j(n_j -\alpha)[2 \psi_0(1-\alpha)\psi_0(n_j -\alpha+1) +\psi_0(n_j-\alpha+1)^2 \right.\nonumber\\
&+&\left.\left.\psi_1(n_j-\alpha+1)+\psi_0(1-\alpha)^2 +\psi_1(1-\alpha)]+2(1-\alpha)[\psi_0(2-\alpha)^2 +\psi_1(2-\alpha) \vphantom{\sum_j}]\right] \right\}\nonumber\\
%\end{eqnarray}
%\begin{eqnarray}
%\label{secondopezzo}
&-&\frac{2}{V_{n,k}}\left\{ \sum_j (n_j-\alpha)[2V_{n+1,k}^* \psi_0(n_j-\alpha +1)+V_{n+1,k}(\psi_0 (n_j-\alpha +1)^2 +\psi_1(n_j -\alpha +1)) +V_{n+1,k}^{**}]\right.\nonumber\\
&+&\left. 2V_{n+1, k+1}^* \psi_0(1-\alpha) +V_{n+1,k+1}[\psi_0(1-\alpha)^2 +\psi_1(1-\alpha)] + V_{n+1,k+1}^{**} \vphantom{\sum_j}\right\}.
\end{eqnarray}
}
%&+&\frac{2}{V_{n,k}} [V_{n+1,k+1}^*+V_{n+1,k}^* \sum_j (n_j-\alpha)+ V_{n+1,k}\sum_j(n_j -\alpha)\psi_0(n_j -\alpha+1)+V_{n+1, k+1}\psi_0(1-\alpha)]
\end{proposition}
\subsection{Examples}
\begin{example}[Two-parameter Poisson-Dirichlet continued] An application of (\ref{post_tsagibbs}) and (\ref{posttsa2}) under two-parameter Poisson-Dirichlet priors yields
\begin{equation}
\label{postmeanAT}
E_{\alpha, \theta, n} (H_m)=\frac{1}{m -1} \left(1 - \left[ \frac{\sum_{j}(n_j -\alpha)_{m}}{(\theta +n)_m}+ \frac{(1 -\alpha)_{m-1}(\theta +k\alpha)}{(\theta +n)_{m}}\right]\right)
\end{equation}
and
{\small
\begin{eqnarray}
\label{postvarAT1}
var_{\alpha, \theta, n}(H_m)&=&\left(\frac{1}{m-1}\right)^2 \times
\left\{ \left[\frac{\sum (n_j -\alpha)_{2m}}{(\theta +n)_{2m}}+ \frac{2 \sum_{i \neq j}(n_j -\alpha)_m(n_i-\alpha)_m}{(\theta +n)_{2m}}\right.\right. \nonumber\\
&+&\frac{(\theta +k\alpha)}{(\theta +n)_{2m}} \left[ (\theta+k\alpha +\alpha)((1-\alpha)_{m-1})^2 +(1-\alpha)_{2m -1}\right]\nonumber\\
&+&\left.\left. \frac{2 \sum (n_j -\alpha)_m (\theta +k\alpha)(1-\alpha)_{m-1}}{(\theta +n)_{2m}}\right]- \left( \frac{\sum_j (n_j -\alpha)_m + (1-\alpha)_{m-1}(\theta +k\alpha)}{(\theta +n)_m}\right)^2\right\}.
\end{eqnarray}
}
For $m=2$ (\ref{postmeanAT}) and (\ref{postvarAT1}) reduce to Simpson's index posterior mean and variance. (See Cerquetti, 2012 for details). Shannon entropy posterior estimation
%relying on well-known stick-breaking representation of the Poisson-Dirichlet size biased atoms$\tilde{P}_j=V_{j}\prod_i (1-V_i)$ for   $\tilde{V}_j \sim Be(1-\alpha, \theta + j\alpha)$, 
arises by an application of (\ref{postSHANgibbs}). For $V_{n,k}=(\theta +\alpha)_{k-1 \uparrow \alpha}/(\theta +1)_{n-1}$ then
\begin{eqnarray}
\label{postmeanSHA}
E_{\alpha, \theta}(H_1|{\bf n})&=&\psi_0(\theta + n +1)\nonumber\\
&-&\frac{1}{(\theta +n)} \left[(\theta +k\alpha) \psi_0(1-\alpha) + \sum_{j=1}^k (n_j -\alpha)\psi_0(n_j -\alpha +1)\right]
\end{eqnarray}
Posterior variance may be obtained through the second posterior moment, by applying (\ref{secpostsha}) for 
$$
\frac{V_{n +\xi m, k+k^*}}{V_{n,k}}= \frac{(\theta +k\alpha)_{k^* \uparrow \alpha}}{(\theta +n)_{\xi m}}
$$
and bearing in mind that
$$
\frac{\partial}{\partial m} \frac{1}{(\theta +n)_{\xi m}}= - \frac{\xi}{ (\theta + n)_{\xi m}} \psi_0 (\theta +n +\xi m)
$$
and
$$
\frac{\partial^2}{\partial^2 m} \frac{1}{(\theta +n)_{\xi m}}= \frac{\xi^2}{(\theta +n)_{\xi m}}[\psi_0 (\theta +n +\xi m)^2 -\psi_1(\theta +n +\xi m)].
$$
It is easy to check that the results agree with Archer {et al.} (2013). Like for prior moments, posterior estimation under Dirichlet priors and normalized Stable priors arises specializing the previous formulas respectively for $\alpha=0$ and $\theta=0$. 
\end{example}
\begin{example} [Gnedin-Fisher priors - continued] Under two-parameter Gnedin Fisher model $(\psi,\gamma)$ Tsallis posterior mean corresponds to
\begin{eqnarray}
\label{postmeanGF2}
E_{\psi, \gamma, {\bf n}}(H_m)=&& \frac{1}{m-1}\times \left[1 - \frac{(\gamma+n-k)_m}{(\psi +n)_m(\gamma+\psi+n)_m}\sum_{j=1}^k(n_j+1)_m \right.\nonumber\\
&-&\left.\frac{(\gamma+n-k)_{m-1}(k-\psi)(k-\gamma+\psi)}{(\psi+m)_m(\gamma-\psi+n)_m}(2)_{m-1}\right].
\end{eqnarray}
Posterior variance follows by an application of (\ref{posttsa2}) for 
$$
\frac{V_{n+\xi m, k+i}}{V_{n,k}}= \frac{(\gamma +n -k)_{\xi m-i} (k-\psi)_i(k-\gamma+\psi)_i}{(\psi+n)_{\xi m}(\gamma+\psi +n)_{\xi m}}.
$$
For $m=2$ previous formulas yield Simpson's index posterior first and second moments.  As for Shannon entropy estimation, posterior mean is easily obtained from (\ref{postSHANgibbs})
\begin{eqnarray}
\label{postshaGF}
E_{\psi, \gamma, {\bf n}}(H_1)&=&\psi_0(\gamma -\psi+n+1)+\psi_0(\psi+n+1)-\psi_0(\gamma+n-k)\nonumber\\ 
&-&\frac{(\gamma +n-k) \sum_j (n_j+1)\psi_0(n_j+2)}{(\psi+n)(\gamma -\psi+n)}\nonumber\\
&-&\frac{(n+k)}{(\psi+n)(\gamma -\psi+n)} -\psi_0(2)\frac{(k-\psi)(k-\gamma +\psi)}{(\psi+n)(\gamma -\psi +n)}.
\end{eqnarray}

Posterior variance may be obtained by (\ref{secpostsha}) and bearing in mind that 
\begin{eqnarray}
\label{ultimafor}
\frac{\partial}{\partial m}\frac{(\gamma +n-k)_{\xi m -i}}{(\psi+n)_{\xi m}(\gamma +\psi +n)_{\xi m}}&=& 
\frac{(\gamma +n-k)_{\xi m -i}}{(\psi+n)_{\xi m}(\gamma +\psi +n)_{\xi m}} \nonumber\\
&\times & \xi [\psi_0 (\gamma +n-k+ \xi m -i)- \psi_0(\psi +\gamma +n +\xi m)- \psi_0(\psi +n+\xi m)]\nonumber
\end{eqnarray}
and 
%\end{document}
\begin{eqnarray}
&&\frac{\partial^2}{\partial^2 m} \frac{(\gamma +n-k)_{\xi m -i}}{(\psi+n)_{\xi m}(\gamma +\psi +n)_{\xi m}}= \xi^2 \frac{(\gamma +n-k)_{\xi m -i}}{(\psi+n)_{\xi m}(\gamma +\psi +n)_{\xi m}}\nonumber\\
&\times & [-2\psi_0(A) (\psi_0(C)+\psi_0(B))+ \psi_0(A)^2 +\psi_1(A) +\psi_0(C)^2 +2\psi_0(B)\psi_0(C)\nonumber\\
&-&\psi_1(C)+\psi_0(B)^2 -\psi_1(B)],
\end{eqnarray}
for $A=(\gamma+n-k +\xi m -i)$, $B=(\psi+n+\xi m)$ and $C=(\gamma+\psi+n+\xi m)$.
\end{example}

\section{Illustration}
\subsection{A real data application}

The high level of generality of the results proposed in the previous sections does not allow a complete direct evaluation of our technique with respect to single alternative frequentist and Bayesian procedures already available in the literature.  Nevertheless a sample illustration of the implementation of our method and a comparison with the results provided by existing procedures can be conducted for a single index in the Tsallis class and a specific prior choice.  Here we apply our Bayesian nonparametric estimation of Shannon index to a dataset on tropical foliage insects from sweep samples taken in 25 sites in Costa Rica and the Caribbean Islands (Janzen, 1973) already considered in Chao and Shen (2003). The dataset consists of the frequency counts for beetles collected in day time from the site referred to as "Osa primary-hill, dry season, 1967". Frequency counts, or sampling formulas, provide an alternative codification of a realization of a random partition in terms of the vector of the number of blocks of the same size. For $m_j$ the number of species represented $j$ times in the beetles sample, the observed dataset is given by $m_1=59, m_2=9, m_3=3, m_4=2, m_5=2, m_6=2, m_{11}=1$, with a total of $k=78$ different species seen and $n=127$ total observations. Most of the species observed have only one, two or three individuals represented in the sample, with few abundant species.  
%Table xx provide  Bayesian point estimation under quadratic loss function (posterior mean) and interval estimation by approximate credibiliy intervals. ????? cosa mostarno questi risultati?  which show that the posterior law of $H_m$ is enough concentrated around the posterior mean.   ??? 

\begin{table}[ht]
\def~{\hphantom{0}}
%\centerin
{\it Point and interval estimates of Shannon's index for beetles dataset}.

%\resizebox{14.5cm}{!}{
\begin{tabular}{rrrr}
  %\hline
Method && Point estimation  & Interval estimation  \\\\
  %\hline
  %$\gamma$\\\\

Maximum likelihood (ML) && 4$\cdot$08 &  (3$\cdot$94, 4$\cdot$22) \\ 
  Bias-corrected ML && 5$\cdot$11 & (4$\cdot$36, 5$\cdot$85) \\ 
  Jackknife && 4$\cdot$62 & (4$\cdot$40, 4$\cdot$84)  \\ 
  Non-parametric ML && 4$\cdot$07  \\ 
  Chao-Shen && 4$\cdot$70 & (4$\cdot$29, 5$\cdot$11) \\\\
BNP - Gnedin-Fisher ($\gamma=$ 0$\cdot$05) && 4$\cdot$859  & (4$\cdot$590, 5$\cdot$136) \\ 
BNP - Gnedin-Fisher ($\gamma=$ 0$\cdot$1)&& 4$\cdot$859  & (4$\cdot$590, 5$\cdot$141) \\ 
BNP - Gnedin-Fisher ($\gamma=$ 0$\cdot$2) && 4$\cdot$856  & (4$\cdot$586, 5$\cdot$126) \\ 
BNP - Gnedin-Fisher ($\gamma=$ 0$\cdot$3) && 4$\cdot$856  & (4$\cdot$589, 5$\cdot$130) \\ \\

   %\hline
\end{tabular}
%}
\label{confrSHANNON}
\end{table}

Table~\ref{confrSHANNON} is an elaboration of Table 6 in Chao and Shen (2003). Those authors propose a nonparametric method combining the Horvitz-Thompson estimator adjusted for unequal probability sampling scheme  and the concept of sample coverage to adjust for the presence of unseen species. They show by simulations that their estimator is preferable to previous frequentist estimators and performs well when a large fraction of the species is missing in the sample.  We add the four estimates obtained applying our method under Gnedin-Fisher priors for the four different values of the $\gamma$ parameter considered in Section 2.3. Bayesian nonparametric estimators (posterior means) are obtained applying formula (\ref{postmeanSHA}), while highest posterior density intervals are derived by the posterior sampled values, by means of a specific R function. With respect to the Horvitz--Thompson estimator corrected for the number of unseen species in Chao and Shen (2003) our technique yields higher point estimates with narrower intervals, thus providing a better account of the effect of missing species with an increased precision. The information in this particular dataset greatly overcomes the information contained in the chosen prior leading to robust conclusions on the diversity of the population, independently of the choice of the prior $\gamma$.  We stress here that, unlike the method proposed by Chao and Shen,  the nonparametric Bayesian approach we are proposing naturally takes into account the presence of unseen species. In fact by construction the prior is placed on a theoretically infinite space of relative abundances for $\alpha \in [0,1)$ and on a random number of finitely many species for $\alpha <0$. Then the multiplicities of the first $k$ species in the observed  sample update both the relative abundances of the seen and of the unseen species. 
\begin{figure}
\centering
\resizebox{14.5cm}{!}{
\begin{tabular}{rrrr}
\includegraphics{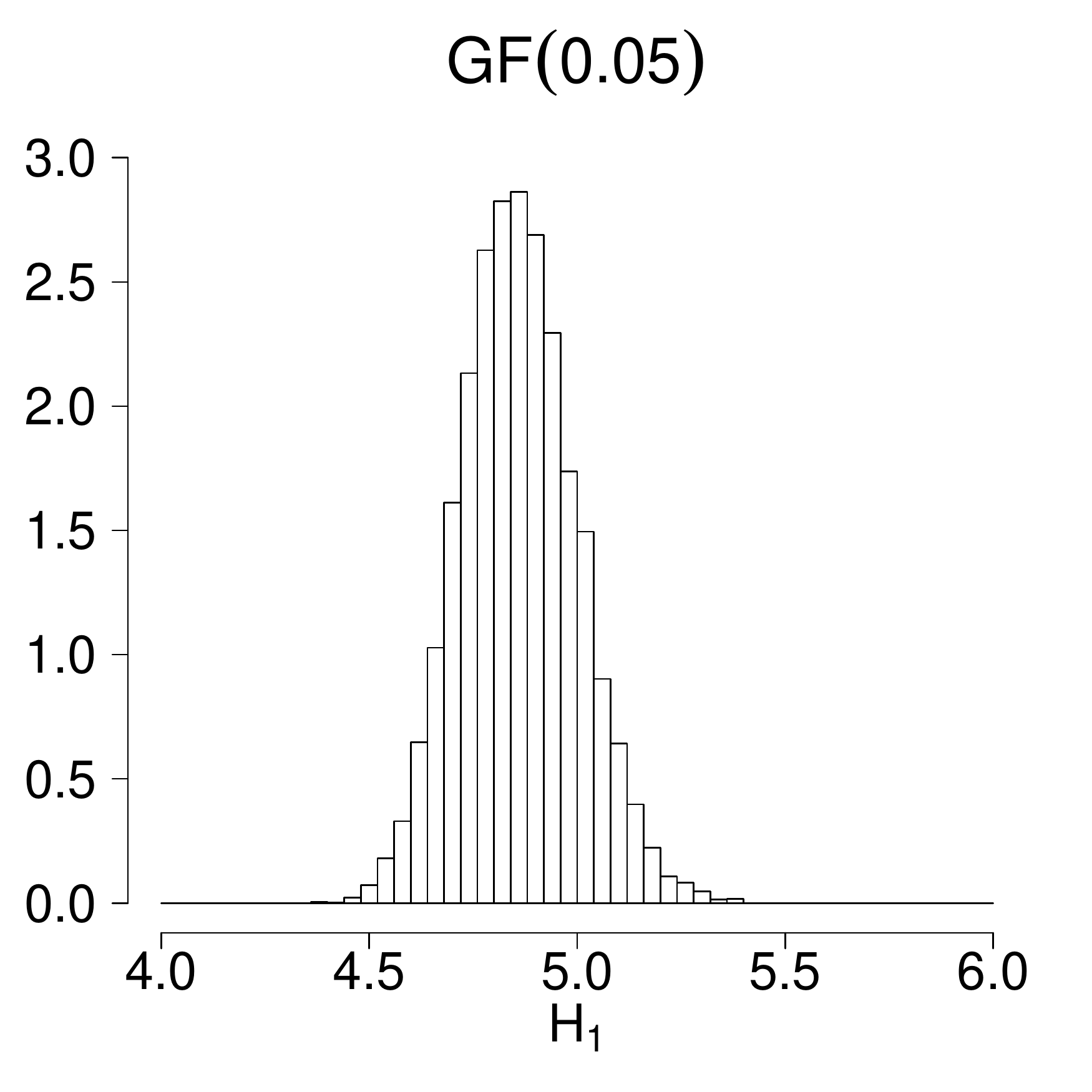}&
\includegraphics{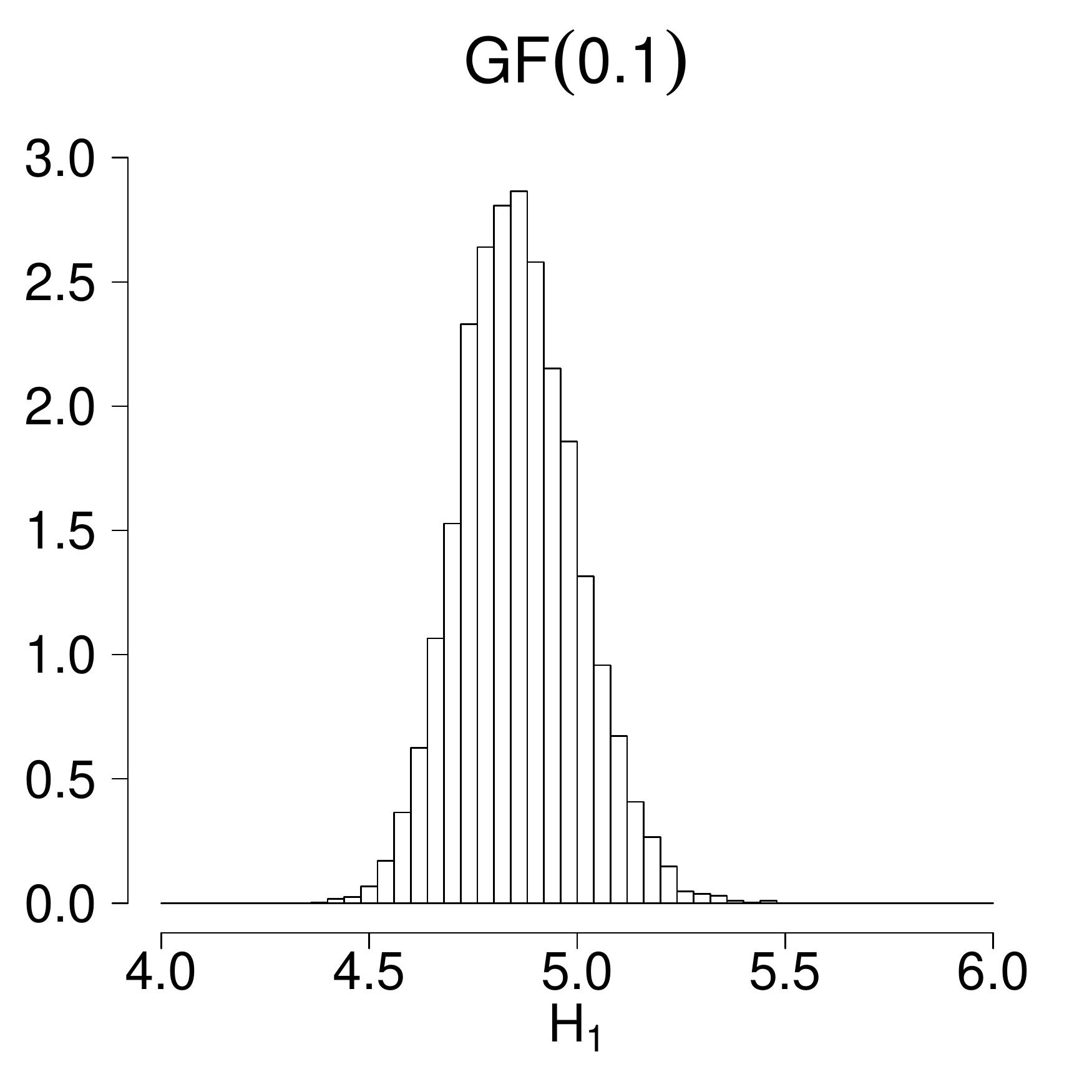}&
\includegraphics{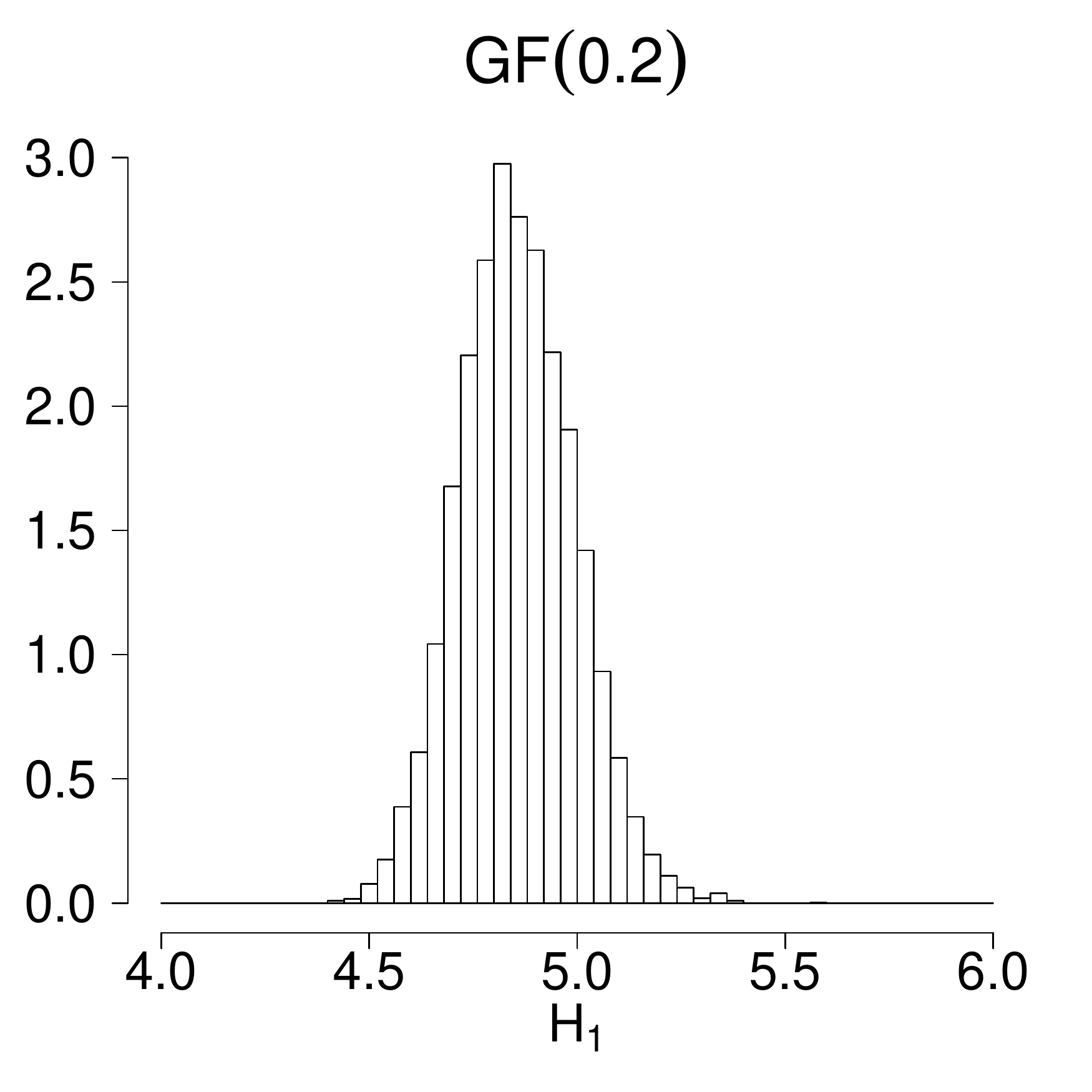}&
\includegraphics{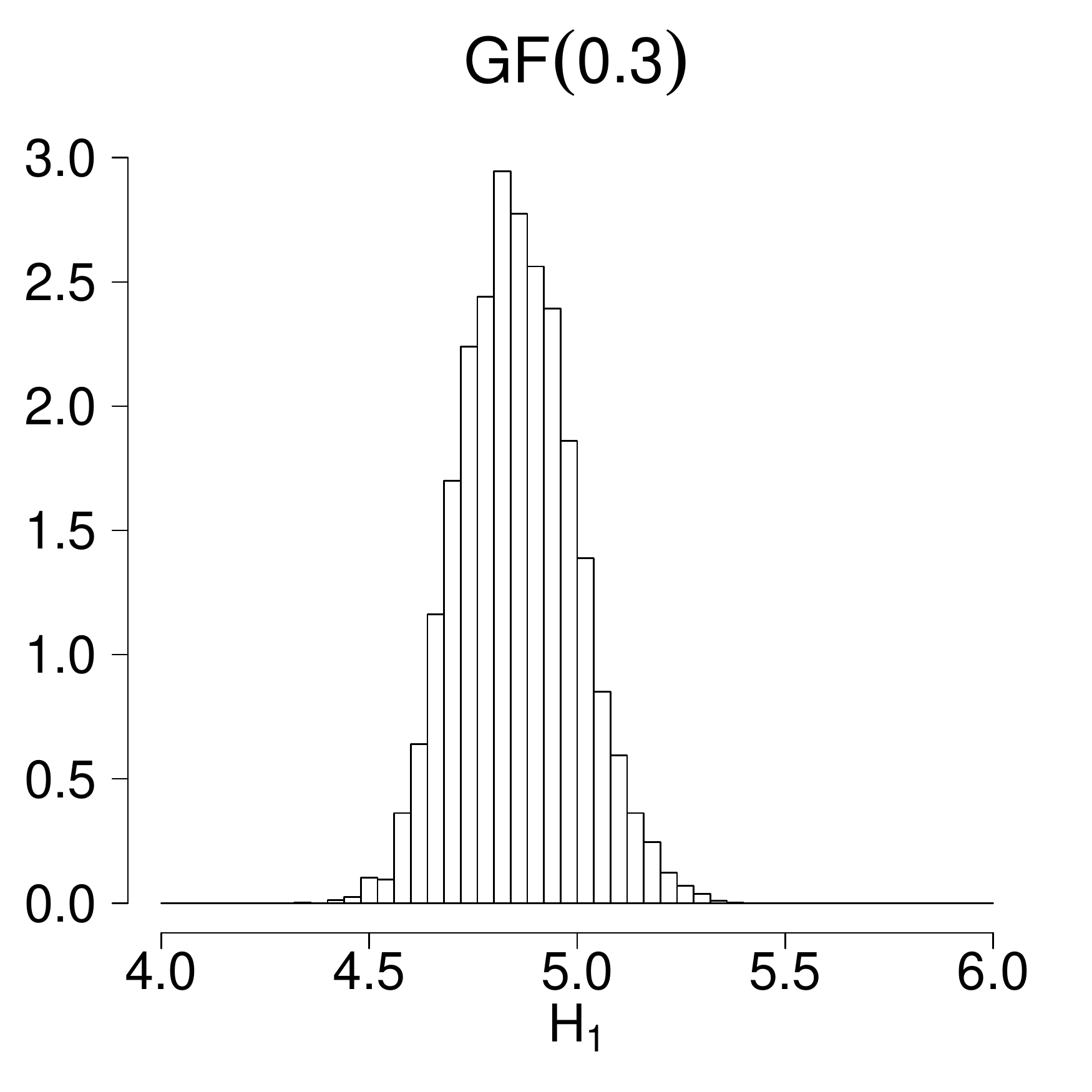}
\end{tabular}
}
\caption{Sampled posterior distributions of $H_1$ under one parameter $\gamma$ Gnedin--Fisher priors for beetles in day time dataset (Janzen, 1973) }
\label{histpriorpostGNEDINFISHER}
\end{figure}
Figure~\ref{histpriorpostGNEDINFISHER} shows the whole posterior distributions under Gnedin-Fisher priors of Shannon index  as obtained by simulations for the beetles dataset. The bell-like shape with a very low level of asymmetry suggests that first and second posterior moments, therefore posterior mean and highest posterior density intervals, are enough to summarize posterior inference on Shannon diversity index in this case.

\section*{Acknowledgement} The authors wishes to thank Leopoldo Catania for his kind assistance in the development of the R code used in the paper, Mauro Bernardi for providing the R function to obtain highest posterior density intervals in Table 2 and Stephan Poppe for introducing her to the notion of Tsallis generalized index. 

%\section*{Supplementary material}
%\label{SM}
%Supplementary material available at \Bka\ online provides the R code used for simulations and results in Section 2 and Section 4.
%\section{Appendix}
\appendix
%\app\endixone
\section*{Appendix }
\subsection{Proofs of Section 2.}
\begin{proof}[of Theorem 1] Let $(P_i^{\downarrow})_{i \geq 1}$ be the sequence of ranked atoms of a random discrete distribution, and $\tilde{P}_j$ the random size of the $j$th atom discovered in the process of random sampling, or equivalently the asymptotic frequency of the $j$th class when the blocks of the partition generated are put in order of their least element.  Now for the random variable 
$$
S_m:=\sum_{i=1}^\infty P_i^m= \sum_{j=1}^\infty \tilde{P}_j^m,
$$
where it is still assumed that $S_1=1$ almost surely, Pitman (2003) provides the following general expression for the $\xi$-th  moment 
\begin{equation}
\label{priormom}
E[(S_m)^\xi]= \sum_{j=1}^\xi \frac{1}{j!} \sum_{\xi_1, \dots, \xi_j} \frac{\xi!}{\xi_1! \cdots \xi_j!} p(m\xi_1, \dots, m\xi_j), 
\end{equation}
where the second sum is over all sequences of $j$ positive integers $(\xi_1, \dots,\xi_j)$ with $\xi_1+\dots+\xi_j=\xi$. For $p(n_1, \dots, n_k)=V_{n,k} \prod_{j=1}^k (1 -\alpha)_{n_j-1}$ (\ref{priormom}) specializes as
\begin{equation}
\label{momsim}
E_{\alpha, V_{n,k}}[(S_m)^\xi]=\sum_{j=1}^\xi \frac{1}{j!} V_{m\xi, j} \sum_{\xi_1, \dots, \xi_j} \frac{\xi!}{\xi_1! \cdots, \xi_j!} \prod_{i=1}^j (1 -\alpha)_{m\xi_i-1}.
\end{equation}
This implies the EPPF induced by sampling from a random discrete distribution directly determines the positive integers moments of the power sums $S_m$, hence the distribution of $S_m$ for each $m$. Explicit first, second and third moments of $S_m$ follows from (\ref{momsim}) for $\xi=1$, $\xi=2$ and $\xi=3$ hence 
\begin{equation}
\label{priormeanSm}
E(S_m)= V_{m,1} (1 -\alpha)_{m-1},
\end{equation}
\begin{equation}
\label{priorsecmomSm}
E[(S_m)^2]= V_{2m, 1}(1-\alpha)_{2m-1}+V_{2m,2}[(1-\alpha)_{m-1}]^2
\end{equation}
and
\begin{equation}
\label{priorthirdmomSm}
E[(S_m)^3]=V_{3m,1}(1-\alpha)_{3m-1}+3V_{3m,2}[(1-\alpha)_{m-1}(1-\alpha)_{2m-1}]+ V_{3m,3}[(1-\alpha)_{m-1}]^3
\end{equation}
and (\ref{mean_tsapr}), (\ref{sec_tsapr}) and (\ref{mom_ter}) easily follow.
\end{proof}
\begin{proof} [of Theorem 2] By definition
$H_1=\lim_{m \rightarrow 1} H_{m}$. By applying H\^opital rule to (\ref{mean_tsapr}) and recalling that $\frac{d}{dx}(\Gamma(x))=\Gamma(x)\psi_0(x)$ (\ref{SHApr_mean}) easily follows. As for the second moment 
$$
\lim_{m \rightarrow 1} E[(H_m)^2]= \lim_{m \rightarrow 1} \frac{1}{(m-1)^2}E\left[(1-S_m)^2\right],
$$
and a repeated application of H\^opital rule yields
$$
\lim_{m \rightarrow 1}\frac{1}{(m-1)^2} [1 +E(S_m)^2 -2E(S_m)]=\lim_{m \rightarrow 1}\frac{1}{2} \left[\frac{\partial^2}{\partial m^2}E[(S_m)^2] -2\frac{\partial^2}{\partial m^2}E(S_m)\right].
$$
Now by  (\ref{priormeanSm}) and (\ref{priorsecmomSm}) 
$$
\frac{\partial^2}{\partial m^2}E[(S_m)^2]= \frac{\partial^2}{\partial m^2}[V_{2m,1}(1-\alpha)_{2m,1}]+ \frac{\partial^2}{\partial m^2} [V_{2m,2}[(1-\alpha)_{m-1}]^2]
$$
and 
$$
\frac{\partial^2}{\partial m^2}E(S_m)= \frac{\partial^2}{\partial m^2} [V_{m.1}(1-\alpha)_{m-1}].
$$
Recalling the definition of {\it trigamma function} $\psi_1(x)=\frac{d}{dx}\psi_0(x)$ then
%$$
%\frac{\partial\frac{\partial \frac{V(m)
%$$
$$
\frac{\partial^2}{\partial m^2}[V_{2m,1}(1-\alpha)_{2m,1}]=
$$ 
$$=\frac{\Gamma(2m-\alpha)}{\Gamma(1-\alpha)} [4 \psi_0(2m-\alpha)\frac{\partial}{\partial m} V_{2m,1}+ 4V_{2m,1}\psi_0(2m-\alpha)^2 + 4 V_{2m,1} \psi_1(2m-\alpha)+ \frac{\partial^2}{\partial m^2}V_{2m,1}],
$$
%$$
%\frac{\partial^2}{\partial m^2}[V_{2m,1}(1-\alpha)_{2m,1}]=\frac{\Gamma (2 m-a) \left(4 \psi ^{(0)}(2 m-a)
%V'(m)+4 V(m) \left(\psi ^{(0)}(2 m-a)^2+\psi ^{(1)}(2 m-a)\right)+V''(m)\right)}{\Gamma (1-a)}
%$$
$$
\frac{\partial^2}{\partial m^2} [V_{2m,2}[(1-\alpha)_{m-1}]^2]= 
$$
$$\frac{[\Gamma(m-\alpha)]^2}{[\Gamma(1-\alpha)]^2}[ 4 \psi_0(m-\alpha)\frac{\partial }{\partial m} V_{2m,2}+ 4 V_{2m,2}\psi_0 (m-\alpha)^2 +2 V_{2m,2}\psi_1(m-\alpha) + \frac{\partial^2}{\partial m^2}V_{2m,2}].
$$
and
$$
 \frac{\partial^2}{\partial m^2} [V_{m.1}(1-\alpha)_{m-1}]= 
$$
$$\frac{\Gamma(m-\alpha)}{\Gamma(1-\alpha)} [2 \psi_0(m-\alpha)\frac{\partial}{\partial m} V_{m,1} + V_{m,1}\psi_0(m-\alpha)^2]+V_{m,1}\psi_1(m-\alpha)+ \frac{\partial^2}{\partial m^2} V_{m,1}.
$$
%$V^*_{2,2}=\lim_{m \rightarrow 1} \frac{\partial}{\partial m} V_{2m,2} $
%$V^{**}_{2,2}=\lim_{m \rightarrow 1} \frac{\partial^2}{\partial m^2} V_{2m,2} $
By taking the limit for $m \rightarrow 1$ (\ref{priSHAgibbs2}) follows.
\end{proof}

\subsection{Proofs of Section 3.}
\begin{proof} [of Theorem 3] 
 Let $(X_i)_{i \geq 1}$ be an exchangeable random sequence of species observations driven by a general $(\alpha, V)$ Gnedin-Pitman prior on the unknown relative abundances and ${\bf n}=(n_1, \dots, n_k)$ the multiplicities of the first $k$ species observed in a sample of size $n$, then for $\xi \geq 1$ posterior moments of $S_m=\sum_{j=1}^{\infty} P_j^m$ can be decomposed as follows
\begin{eqnarray}
\label{postmom}
&&E_{P|{\bf n}}[(\sum_{j=1}^\infty \tilde{P}_j^m)^\xi]=\nonumber\\
&=&E_{P| \bf{n}}[(\sum_{j=1}^k \tilde{P}_j^m)^\xi]+ \sum_{l=1}^{\xi-1}{\xi \choose l} E_{P| \bf{n}}[(\sum_{j=1}^k \tilde{P}_j^m)^{\xi-l}(\sum_{j=k+1}^{\infty} \tilde{P}_j^m)^l]+E_{P| \bf{n}}[(\sum_{j=k+1}^\infty \tilde{P}_j^m)^\xi]
\end{eqnarray} 
By easy combinatorics and telescoping product from the one-step prediction rules under (\ref{EPPFgibbs}) the conditional probability of any {\it particular} partition of the set $[n+v]-[n]$ in $k^*$ new blocks of size $s_i \geq 1$, $\sum_{i=1}^{k^*} s_i=s$, $s \leq v$, with allocation in $k$ old blocks of  $m_j \geq 0$, $\sum_{j=1}^{k} m_j=v-s$ integers, corresponds to
\begin{equation}
\label{oldenew}
%\label{gibbsoldnew}
p_{\bf m}^ {\bf s}({\bf n})=\frac{V_{n+v,k+k^*}}{V_{n,k}}\prod_{j=1}^k (n_j-\alpha)_{m_j}\prod_{i=1}^{k^*}(1-\alpha)_{s_i-1}, 
\end{equation}
for ${\bf n}=(n_1, \dots, n_k)$, ${\bf s}=(s_1, \dots, s_{k^*})$ and ${\bf m}=(m_1, \dots, m_k)$.  For $k^*=0$ then 
\begin{equation}
\label{old}
%\label{gibbsoldnew}
p_{\bf m}({\bf n})=\frac{V_{n+v,k}}{V_{n,k}}\prod_{j=1}^k (n_j-\alpha)_{m_j} 
\end{equation}
and for $v=\sum_i s_i$ then
\begin{equation}
\label{new}
%\label{gibbsoldnew}
p^{\bf s}({\bf n})=\frac{V_{n+v,k+k^*}}{V_{n,k}}\prod_{i=1}^{k^*}(1-\alpha)_{s_i-1}. 
\end{equation}
For $v=m\xi$ and $m_j=m\xi_j$  from (\ref{old}) and multinomial formula,
\begin{equation}
\label{primikappa}
E_{P| \bf{n}}[(\sum_{j=1}^k \tilde{P}_j^m)^\xi]=\sum_{(\xi_1, \dots, \xi_k)} \frac{\xi!}{\prod_j \xi_j!} \frac{V_{n+\xi m, k}}{V_{n,k}} \prod_{j=1}^k (n_j-\alpha)_{m\xi_j}, 
\end{equation}
where the sum is over the space of non-negative integers $(\xi_1, \dots, \xi_k)$ with sum $\xi$. For $v=m\xi$ and $s_i=m\xi_j$ then, from an application of (\ref{momsim}) to the posterior partition probability function (\ref{new}), 
\begin{equation}
\label{resto}
E_{P| \bf{n}}[(\sum_{j=k+1}^\infty \tilde{P}_i^m)^\xi]= \sum_{k^*=1}^\xi \frac{1}{k^*!}\frac{V_{n+m\xi, k+k^*}}{V_{n,k}}\sum_{(z_1, \dots, z_{k^*})} \frac{\xi!}{\prod_{i} z_i!} \prod_{i=1}^{k^*}(1-\alpha)_{mz_i-1}.
\end{equation}
and analogously from  (\ref{oldenew})
\begin{eqnarray}
\label{misti}
&&E_{P| \bf{n}}[(\sum_{j=1}^k \tilde{P}_j^m)^{\xi-l}(\sum_{j=k+1}^{\infty} \tilde{P}_j^m)^l]=\nonumber\\
&=&\sum_{k^*= 1}^ l
\frac{1}{k^*!} \frac{V_{n+m \xi, k+k^*}}{V_{n,k}} \sum_{(\xi_1, \dots, \xi_k)} \frac{\xi -l !}{\prod_j \xi_j!} \prod_{j=1}^k (n_j-\alpha)_{m \xi_j} \sum_{(z_1, \dots, z_{k^*})} \frac{l!}{\prod_{i} z_i!} \prod_{i=1}^{k^*} (1-\alpha)_{mz_i -1},
\end{eqnarray}
where the second sum is over the space of positive integers  $(z_1, \dots, z_{k^*})$ with sum $l$. Suitably applying (\ref{primikappa}), (\ref{resto}) and (\ref{misti}), then (\ref{post_tsagibbs}), (\ref{posttsa2}) and (\ref{postSm3}) follow. 
%An application of (\ref{priormom}) to (\ref{oldandnew})(\ref{old}) (with a fixed number of blocks) and to (\ref{new}) yields the result.
\end{proof}

\section*{References}
\newcommand{\bibu}{\item \hskip-0.4cm}
\begin{list}{\ }{\setlength\leftmargin{0.5cm}}

\small{
\bibu \textsc{Archer, E., Park, M.I. and Pillow, J.} (2013) Bayesian entropy estimation for countable discrete distributions.  	{\it arXiv:1302.0328 [cs.IT]}

\bibu \textsc{Butturi-Gomes, D., Junior, P.M., Giacomini, H.C. and Junior P.D.M.} (2014) Computer intensive methods for controlling bias in generalized species diversity index. {\it Ecological Indicators} {\bf 37}, 90--98.

\bibu \textsc{Cerquetti, A.} (2011) On some Bayesian nonparametric estimators for species richness under two-parameter Poisson-Dirichlet priors. {\it arXiv:1002.0535 [math.PR]}

\bibu \textsc {Cerquetti, A.} (2011b) Reparametrizing the two-parameter Gnedin-Fisher partition model in a Bayesian perspective. {\it Proceedings of the ISI 58th World Statistical Congress}, 2011, Dublin. 4678--4683.

\bibu \textsc{Cerquetti, A.} (2012) Bayesian nonparametric estimation of Simpson's evenness index under $\alpha$-Gibbs priors. {\it arXiv:1203.1666}

\bibu \textsc{Cerquetti, A.} (2013) Marginals of multivariate Gibbs distributions with applications in Bayesian species sampling {\it Elect. J. Stat.},7,  697--716. 

\bibu \textsc{Cerquetti, A.} (2013b) Yet another application of marginals of multivariate Gibbs distributions. {\it arXiv:1312.5789 [stat.ME]}

%\bibu \textsc{Cerquetti, A.} (2013c) A note on a Bayesian estimator of the discovery probability. {\it arXiv:1304.1030 [math.ST]}

\bibu \textsc {Chao, A. Shen, T.} (2003) Nonparametric estimation of Shannon's index of diversity when there are unseen species in sample. {\it Envirom. Ecolog. Statistics}, {\bf 10}, 429--443.

\bibu \textsc{Ewens, W. J.} (1972)  The sampling theory of selectively neutral alleles. {\it Theor. Pop. Biol.}, 3, 87-112.

\bibu \textsc{Favaro, S., Lijoi, A., Mena, R.H. and Pr\"unster, I.} (2009) Bayesian non-parametric inference for species variety with a two-parameter Poisson-Dirichlet process prior. {\it JRSS-B}, {\bf 71}, 993-1008.

\bibu \textsc {Favaro, S., Lijoi, A. and Pr\"unster, I.} (2012) On the stick breaking representation of normalized inverse Gaussian priors. {\it Biometrika}, {\bf 99}, 663--674.

\bibu \textsc{Favaro, S. Lijoi, A. and Pr\"unster, I.} (2012b) A new estimator of the discovery probability. {\it Biometrics}, {\bf 68}, 1188--1196.

\bibu \textsc{Favaro, S. Lijoi, A. and Pr\"unster, I.} (2013) Conditional formulae for Gibbs-type exchangeable random partitions. {\it Ann. Appl. Probab.}, {\bf 23}, 1721--1754.

\bibu \textsc{Ferguson, T.S.} (1973) A Bayesian analysis of some nonparametric problems. {\it Ann. Statist.}, {\bf 1}, 209--230.

\bibu \textsc{Fisher, R.A., Corbet, A.S. and Williams, C. B.} (1943) The relation between the number of species and the number of individuals in a random sample of an animal population. {\it J. Animal Ecol.} {\bf 12}, 42--58.

\bibu \textsc{Gill, C.A. and Joanes, D. N.} (1979) Bayesian estimation of Shannon's index of diversity. {\it Biometrika}, {\bf 66}, 1, 81-85.

%\bibu \textsc{Ginebra, J.and Puig, X.} (2010) On the measure and the estimation of evenness and diversity. {\it Computational Statistics and Data Analysis}, 54, 2187--2201. 

\bibu \textsc{Gnedin, A.} (2010) A species sampling model with finitely many types.  {\it Electron. Commun. Prob.}, {\bf 15}, 79-88. 

\bibu \textsc{Gnedin, A. and Pitman, J. } (2006) {Exchangeable Gibbs partitions  and Stirling triangles.} {\it Journal of Mathematical Sciences}, {\bf 138}, 3, 5674--5685. 

\bibu \textsc{Holste, D. Grosse, I. and Herzel, H.} (1998)  Bayes' estimators of generalized entropies. {\it J. Phys. A: Math. Gen.} {\bf 31}, 2551--2566.

\bibu \textsc{Janzen, D.H.} (1973) Sweep samples of tropical foliage insects: description of study sites, with data on species abundances and size distributions. {\it Ecology}, {\bf 54}, 659--686.

\bibu \textsc{Jost, L.} (2006) Entropy and diversity. {\it Oikos}, {\bf 113}, 2. 

\bibu \textsc{Keylock, C.J.} (2005) Simpson diversity and the Shannon--Wiener index as special cases of a generalized entropy. {\it Oikos}, {\bf 109}, 203--207.

\bibu \textsc{Lijoi, A., Mena, R.H. and Pr\"uenster, I.} (2005) Hierarchical Mixture Modeling With Normalized Inverse-Gaussian Priors. {\it JASA}, {\bf 100}, 1278--1291.

\bibu \textsc{Lijoi, A., Mena, R.H. and Pr\"unster, I.} (2007) Bayesian nonparametric estimation of the probability of discovering new species.  {\it Biometrika}, {\bf 94}, 769--786.

\bibu \textsc{Lijoi, A., Pr\"unster, I. and Walker, S.G.} (2008) Bayesian nonparametric estimator derived from conditional Gibbs structures. {\it Ann. Appl. Probab.}, {\bf 18}, 1519--1547.

%\bibu \textsc {Lloyd, M and Ghelardi, R. J.} (1964) A table for calculating the equitability component of species diversity. {\it J. Animal Ecol.} 33, 217--225.

\bibu \textsc{Mendes, R.S, Evangelista, L. R., Thomaz, S. M. Agostinho, A. A, Gomes, L.C} (2008) A unified index to measure ecological diversity and species rarity. {\it Ecography}, 31, 4, 450--456.

\bibu \textsc{Martins, A.F.T., Smith, N.A, Xing, E.P., Aguiar, P.M.Q. Figuereido, M.A.T.} (2009) Nonextensive information theoretic kernels on measures. {\it J. of Mach. Learning}, {\bf 10}, 935-975.

\bibu \textsc{Nemenman, I., Bialek, W. and De Ruyter van Steveninck, R.} (2004) Entropy and information in neural spike trains: Progress on the sampling problem. {\it Physical Review E}, {\bf 69}, 056111. 

\bibu \textsc{Nemenman, I., Shafee, F. and Bialek, W.} (2002) Entropy and inference, revisited . {\it Adv, Neur, Inf. Proc. Sys.}, 14.

\bibu \textsc{Patil, G. P. and Taillie, C.} (1982) Diversity as a concept and its measurement. {\it JASA}, {\bf 77}, 548--561.

\bibu \textsc{Pielou, E.C.} (1975) {\it Ecological Diversity} New York: Wiley.

\bibu \textsc{Pitman, J.} (1995) Exchangeable and partially exchangeable random partitions. {\it Probab. Th. Rel. Fields}, {\bf 102}: 145-158.

\bibu \textsc{Pitman, J.} (1996) Some developments of the Blackwell-MacQueen urn scheme. In T.S. Ferguson, Shapley L.S., and MacQueen J.B., editors, {\it Statistics, Probability and Game Theory}, volume 30 of {\it IMS Lecture Notes-Monograph Series}, pages 245--267. Institute of Mathematical Statistics, Hayward, CA.

\bibu \textsc{Pitman, J.} (2003) {Poisson-Kingman partitions}. In D.R. Goldstein, editor, {\it Science and Statistics: A Festschrift for Terry Speed}, volume 40 of Lecture Notes-Monograph Series, pages 1--34. Institute of Mathematical Statistics, Hayward, California.

\bibu \textsc{Pitman, J.} (2006) {\it Combinatorial Stochastic Processes}. Ecole d'Et\'e de Probabilit\'e de Saint-Flour XXXII - 2002. Lecture Notes in Mathematics N. 1875, Springer.

\bibu \textsc{Pitman, J. and Yor, M.} (1997) The two-parameter Poisson-Dirichlet distribution derived from a stable subordinator. {\it Ann. Probab.}, {\bf 25}, 855--900.

\bibu \textsc{Shannon, C.E.} (1948) A mathematical theory of communication. {\it Bell System Technical Journal}, 27, 379--423.

\bibu \textsc {Simpson, E.H.} (1949) Measurement of diversity. {\it Nature} {\bf 163}, 688

\bibu \textsc {Tsallis, C.} (1988) Possible generalization of Boltzmann--Gibbs statistics. {\it Journal of Statistical Physics}, {\bf 52}, 479--487.

\bibu \textsc{Tsallis, C.}  (2009) {\it Introduction to Nonextensive Statistical Mechanics: approaching a complex world}. New York: Springer. 

\bibu \textsc{Vila, M., Bardera, A., Feixas, M. and Sbert, M.} (2011) Tsallis mutual information for document classification. {\it Entropy}, {\bf 13}, 1694--1707.

%\bibu \textsc{Wolpert, D.H. and DeDeo, S.} (2013) Estimating functions of distributions defined over spaces of unknown size. {\it Entropy}, 15, 4668--4699.

\bibu \textsc{Wolpert, D.H.,  Wolf, D. R.} (1995) Estimating functions of probability distributions from a finite set of samples. {\it Physical Review E}, {\bf 52}, 6.

\bibu \textsc{Xekalaki, E.} (1983) Infinite divisibility, completeness and regression properties of the univariate generalized Waring distribution. {\it Ann. Inst. Statist. Math.}, {\bf 35}, 279-289.

\bibu \textsc{Zhang, D., Jia, X., Ding, H, Ye, D. and Thakor, N.V.} (2010) Application of Tsallis entropy to EEG: quantifying the presence of burst suppression after asphyxial cardiac arrest in rats. {\it IEEE Trans Biomed Eng.} {\bf 57} (4) 867--874.

}
\end{list}

\end{document}